\documentclass{amsart}

\RequirePackage{xcolor}
\definecolor{mg}{RGB}{34, 139, 34}

\RequirePackage{float}
\RequirePackage[margin=1in]{geometry}
\RequirePackage{graphicx}
\RequirePackage{amsmath, amsthm, amssymb, amsfonts,slashed}
\usepackage{aliascnt}
\RequirePackage{xifthen}
\RequirePackage{indentfirst}
\RequirePackage{microtype}

\newenvironment{Eq}[1][]{\begin{equation}\ifthenelse{\equal{#1}{}}{}{\tag{#1}}\begin{aligned}}{\end{aligned}\end{equation}\ignorespacesafterend}
\newenvironment{Eq*}[1][]{\begin{equation*}\ifthenelse{\equal{#1}{}}{}{\tag{#1}}\begin{aligned}}{\end{aligned}\end{equation*}\ignorespacesafterend}
\numberwithin{equation}{section}
\newtheorem{theorem}{Theorem}[section]

\newaliascnt{proposition}{theorem}
\newtheorem{proposition}[proposition]{Proposition}
\aliascntresetthe{proposition}

\newaliascnt{corollary}{theorem}
\newtheorem{corollary}[corollary]{Corollary}
\aliascntresetthe{corollary}

\newaliascnt{lemma}{theorem}
\newtheorem{lemma}[lemma]{Lemma}
\aliascntresetthe{lemma}

\newaliascnt{claim}{theorem}
\newtheorem{claim}[claim]{Claim}%
\aliascntresetthe{claim}

\theoremstyle{remark}%
\newtheorem{remark}{Remark}[section]
\theoremstyle{definition}

\usepackage[colorlinks=true,linkcolor=red,citecolor=blue,urlcolor=blue]{hyperref}
\usepackage[capitalise,noabbrev]{cleveref}

\DeclareMathOperator{\fd}{d} \renewcommand{\d}{\fd}

\newcommand{\sm}{{\slashed m}}
\newcommand{\sGamma}{{\slashed \Gamma}}

\newcommand{\rR}{{\mathbb{R}}}
\newcommand{\rH}{{\mathcal H}}
\newcommand{\rHb}{{\underline{\mathcal H}}}
\newcommand{\rS}{{\mathbb{S}}}

\newcommand{\rZ}{{\mathbb{Z}}}
\newcommand{\rD}{{\mathcal D}}

\newcommand{\bn}{{\mathbf n}}
\newcommand{\spartial}{{\slashed \partial}}
\newcommand{\hpartial}[1]{\widehat{\partial_{#1}}}

\newcommand{\kl}[1]{\mathopen{}\left#1}
\newcommand{\kr}[1]{\right#1}
\newcommand{\Low}[1]{\limits_{\substack{#1}}}

\crefname{equation}{}{}
\crefname{theorem}{Theorem}{Theorems}
\Crefname{theorem}{Theorem}{Theorems}
\crefname{proposition}{Proposition}{Propositions}
\Crefname{proposition}{Proposition}{Propositions}
\crefname{corollary}{Corollary}{Corollaries}
\Crefname{corollary}{Corollary}{Corollaries}
\crefname{lemma}{Lemma}{Lemmas}
\Crefname{lemma}{Lemma}{Lemmas}
\crefname{claim}{Claim}{Claims}
\Crefname{claim}{Claim}{Claims}
\crefname{hypothesis}{Hypothesis}{Hypotheses}
\Crefname{hypothesis}{Hypothesis}{Hypotheses}
\crefname{remark}{Remark}{Remarks}
\Crefname{remark}{Remark}{Remarks}
\crefname{definition}{Definition}{Definitions}
\Crefname{definition}{Definition}{Definitions}

\title[Wave systems with derivative nonlinearities]
{Long-time existence for wave systems with derivative nonlinearities}

\author{Wei Dai}
\address{School of Mathematical Sciences\\
Zhejiang University of Technology\\
Hangzhou 310023, P. R. China}
\email{daiw23@zjut.edu.cn}

\author{Mengyun Liu}
\address{Department of Mathematics\\
Zhejiang Sci-Tech University\\
Hangzhou 310018, P. R. China}
\email{mengyunliu@zstu.edu.cn}

\author[Chengbo Wang]{Chengbo Wang$^{*}$}
\thanks{$^{*}$Corresponding author.}
\address{School of Mathematical Sciences\\
Zhejiang Normal University\\
Jinhua 321004, P. R. China}
\email{wangcbo@zjnu.edu.cn}

\date{\today}

\begin{document}
\bibliographystyle{plain}

\begin{abstract}
In this paper, we study wave systems with power-type derivative nonlinearities
in space dimensions $n\geq3$, for small, spherically symmetric initial data.
Under suitable assumptions on the nonlinear powers, we establish global
existence and sharp lifespan lower bounds.
To achieve this, we develop weighted energy and local energy estimates
adapted to the null directions, together with corresponding pointwise
decay estimates.
\end{abstract}

\keywords{Glassey conjecture, semilinear wave systems, lifespan,
weighted energy estimates, local energy estimates}

\subjclass[2020]{Primary 35L71; Secondary 35L52, 35B33, 35B45}

\maketitle

\section{Introduction}\label{sec:introduction}
In this paper, we consider the following system of semilinear wave equations
with derivative-type nonlinearities:
\begin{Eq}\label{eq:main-system}
    \begin{cases}
        \square \phi=\phi_{tt}-\Delta\phi=|\psi_t|^{p},\qquad t>0,x\in\rR^n,\\
        \square \psi=\psi_{tt}-\Delta\psi=|\phi_t|^{q},\qquad t>0,x\in\rR^n,\\
        (\phi,\phi_t,\psi,\psi_t)_{t=0}=\varepsilon(f_\phi,g_\phi,f_\psi,g_\psi).
    \end{cases}
\end{Eq}
Here $n\geq3$, $p,q>1$, and $0<\varepsilon\ll1$.
Our main purpose is to establish small-data global existence and sharp
lifespan lower bounds for spherically symmetric solutions, with particular
emphasis on the higher-dimensional problem.

Before presenting our main result, let us briefly recall the relevant history.
For the scalar equation
\begin{Eq*}
\square \phi=|\phi_t|^p,
\end{Eq*}
Glassey
\cite{GlasseyMR1983} conjectured that the critical exponent separating
small-data global existence from finite-time blow-up is
\begin{Eq*}
p_G(n)=1+\frac{2}{n-1}.
\end{Eq*}
In space dimensions $n=2,3$, the conjecture was verified for general initial
data by Hidano--Tsutaya \cite{MR1386769} and Tzvetkov
\cite{MR1637692}, while an earlier radial result in three dimensions was
obtained by Sideris \cite{Sideris1983}. For $n\geq4$, Zhou
\cite{MR1845748} proved finite-time blow-up together with the expected sharp
upper bounds for the lifespan when $1<p\leq p_G(n)$. The existence side in
higher dimensions was later completed under spherical symmetry by
Hidano--Wang--Yokoyama \cite{MR2980460}, who proved small-data global
existence for $p_G(n)<p<\frac{n}{n-2}$ and obtained sharp lower bounds for the lifespan in
the critical and subcritical regimes. Thus, for the scalar equation, the
critical-threshold and lifespan picture is by now essentially understood
in the radial setting.

By comparison, the literature for the corresponding weakly coupled systems is
much less complete. In three space dimensions, consider the system
\begin{Eq}\label{eq:multiple-speed-derivative-system}
\begin{cases}
\phi_{tt}-c_1^2\Delta\phi=|\psi_t|^p,\\
\psi_{tt}-c_2^2\Delta\psi=|\phi_t|^q,
\end{cases}
\qquad c_1,c_2>0,\quad p,q>1.
\end{Eq}
For $c_1=c_2=1$, we may assume $1<p\leq q$ without loss of generality
by interchanging the two equations. Deng \cite{Deng1999} proved
finite-time blow-up when
\begin{Eq*}
q(p-1)\leq2,
\end{Eq*}
whereas Kubo--Kubota--Sunagawa \cite{MR2221121} established small-data global
existence for radially symmetric initial data when
\begin{Eq*}
q(p-1)>2.
\end{Eq*}
Consequently, the sharp critical curve for the three-dimensional radial
equal-speed problem was determined. Systems with distinct propagation
speeds were also considered in the latter paper, as we discuss briefly
below. Some typical
three-dimensional cases were subsequently treated without the radial
assumption by Kubo--Takaki \cite{KT2010}.

For the equal-speed problem in general space dimensions,
Ikeda--Sobajima--Wakasa \cite{I-S-W} proved finite-time blow-up and obtained
lifespan upper bounds governed by the quantity
\begin{Eq}\label{eq:critical-Gamma}
\Gamma(p,q)
:=\max\kl\{
\frac{p+1}{pq-1},
\frac{q+1}{pq-1}
\kr\}
-\frac{n-1}{2},
\end{Eq}
for $n\geq 2$.
Their blow-up result applies in the region $\Gamma(p,q)\geq0$. This left a rather
asymmetric picture: in higher dimensions the blow-up side, including
quantitative lifespan upper bounds, was available, while
matching existence results and lifespan lower bounds had not been
established; even in three dimensions, where
the critical curve had already been identified by
Kubo--Kubota--Sunagawa \cite{MR2221121}, quantitative lifespan lower bounds in the
critical and subcritical regimes were not obtained there.
In the present paper, we address these questions by establishing radial
small-data global existence and lifespan lower bounds for the
equal-speed system, within the range of exponents specified below.

One reason for restricting our attention to spherically symmetric solutions
is the limited regularity of the nonlinearities in the relevant
higher-dimensional range. Even for the scalar Glassey
equation, sharp higher-dimensional existence results near the critical power
are available essentially under spherical symmetry; see
Hidano--Wang--Yokoyama \cite{MR2980460} and Wang
\cite{MR3378835,Wang2015Obstacles}. One basic difficulty is
that the relevant high-dimensional range involves powers below $2$, for which
the derivative nonlinearities have only limited regularity. Closing the
energy estimates requires suitable pointwise control, but the limited
regularity makes it difficult to exploit the angular regularity used in
nonradial arguments. In the radial setting, the required pointwise bounds
can instead be obtained from Hardy- and trace-type inequalities.

We also mention several related developments. For example, the scalar
theory has been developed on asymptotically flat
manifolds and exterior domains \cite{MR3378835,Wang2015Obstacles},
and finer critical nonlinearities involving a modulus of continuity have also
been studied \cite{ChenPalmieri2023,Shao2024Criteria}.
For the three-dimensional system
\eqref{eq:multiple-speed-derivative-system} with $c_1\neq c_2$,
Kubo--Kubota--Sunagawa \cite{MR2221121} obtained radial global existence
under a condition different from that in the equal-speed case, but did not
determine the sharp critical curve. In a recent preprint,
Dai--Liu--Wang \cite{DaiLiuWangMultipleSpeed} determine the sharp
critical curve for the same system in the range $1<p,q<2$ by combining
the global-existence and blow-up regions.

To state our main result, let
\begin{Eq*}
\boldsymbol f:=(f_\phi,g_\phi,f_\psi,g_\psi).
\end{Eq*}
Define
\begin{Eq*}
\|\boldsymbol f\|_{\mathcal X}
:={}&
\sum_{h\in\{\phi,\psi\}}
  \left\|
  \kl<x\kr>^{\frac32}
  \kl(
  \nabla^{\leq 2}f_h,
  \nabla^{\leq 1} g_h
  \kr)
  \right\|_{L^2(\rR^n)},
\end{Eq*}
with $\kl<x\kr>:=\sqrt{1+|x|^2}$.
Here and below, for any family $Z$ of vector fields, we write
$Z^{\leq k}f:=\{Z^\alpha f\}_{|\alpha|\leq k}$.
Let $\mathcal X^{\mathrm{rad}}$ be the completion of
$\kl(C_{c,\mathrm{rad}}^\infty(\rR^n)\kr)^4$
with respect to this norm, where
$C_{c,\mathrm{rad}}^\infty(\rR^n)$ denotes the subspace of
$C_c^\infty(\rR^n)$ consisting of spherically symmetric functions.

\begin{theorem}\label{thm:main}
Let $n\geq3$ and $1<p,q<n/(n-2)$.
Suppose that $\boldsymbol f\in\mathcal X^{\mathrm{rad}}$ satisfies
$\|\boldsymbol f\|_{\mathcal X}\leq1$.
Then there exist constants $0<\varepsilon_0=\varepsilon_0(n,p,q)\ll 1$ and
$c=c(n,p,q)>0$ such that, for every
$0<\varepsilon\leq\varepsilon_0$, the Cauchy problem
\eqref{eq:main-system} admits a unique
spherically symmetric $H^2$ solution, whose maximal lifespan $T_*$ satisfies
\begin{Eq}\label{eq:lifespan-lower-bound}
T_*\geq T_\varepsilon:=
\begin{cases}
\infty,
&\Gamma(p,q)<0,\\
\exp\kl(c\varepsilon^{-(pq-1)}\kr),
&\Gamma(p,q)=0,\quad p\neq q,\\
\exp\kl(c\varepsilon^{-(p-1)}\kr),
&\Gamma(p,q)=0,\quad p=q,\\
c\varepsilon^{-1/\Gamma(p,q)},
&\Gamma(p,q)>0.
\end{cases}
\end{Eq}
\end{theorem}

\begin{remark}
Comparing our lifespan lower bounds with the upper bounds in \cite{I-S-W},
we see that our estimates are sharp in order.
More precisely, \cref{fig:critical-curve} illustrates the resulting picture.
Within the exponent range $1<p,q<n/(n-2)$, marked by the horizontal and
vertical dashed lines, the region above and to the right of the critical
curve $\Gamma(p,q)=0$ admits global small-data radial solutions.
On the critical curve and in the region below and to its left,
\cref{thm:main} gives exponential and polynomial lifespan lower bounds,
respectively, with the sharp dependence on $\varepsilon$ described above.
\begin{figure}[H]
\centering
\includegraphics[width=0.3\textwidth]{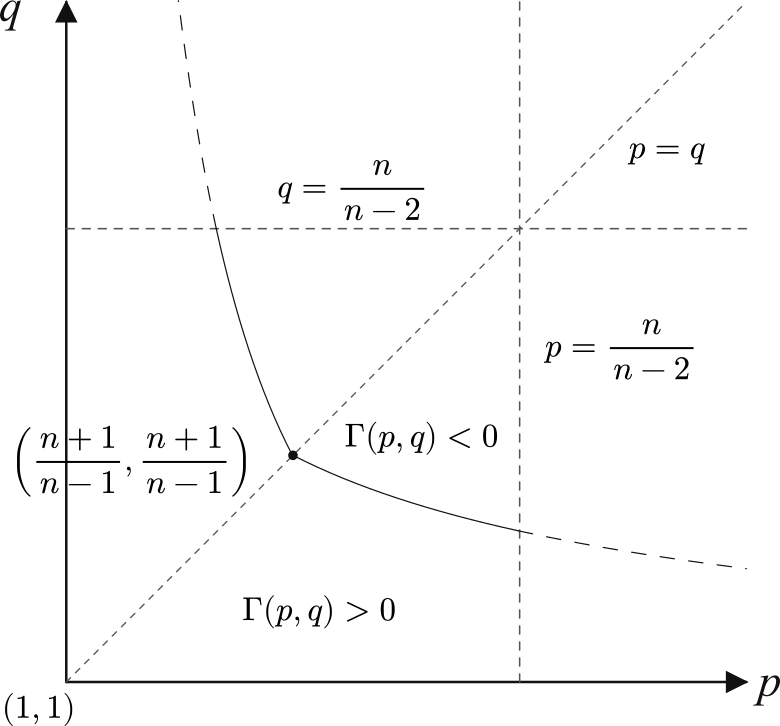}
\caption{The critical curve $\Gamma(p,q)=0$ and the exponent range of
\cref{thm:main}.}
\label{fig:critical-curve}
\end{figure}
\end{remark}

\begin{remark}
Our argument also applies to systems whose
nonlinearities involve radial derivatives such as $\psi_r$ and $\phi_r$, and
yields the same conclusions. We omit this extension to keep the presentation
concise.
\end{remark}

Let us briefly comment on a main difficulty in extending the previous scalar
arguments to the present coupled system. For example, in three space
dimensions and in the regime $1<p<2<q$, the estimates of
Kubo--Kubota--Sunagawa \cite{MR2221121} suggest that a direct adaptation of
the scalar energy framework of Hidano--Wang--Yokoyama \cite{MR2980460}
would allow polynomial growth in time for one of the components. The scalar
estimates alone are insufficient to compensate for this growth, so the
corresponding energy argument does not close directly. In the coupled
problem, one is therefore led to exploit additional time decay from the
other component to compensate for this growth. The purely spatial weights
used in the scalar framework do not by themselves provide the required time
decay. To overcome this difficulty, we use weights adapted to the null
directions together with Hardy-type and trace-type inequalities along the
null cones. This yields the additional time decay needed to close the
coupled estimates.

A main ingredient of the proof is the family of weighted energy and local
energy estimates developed here, adapted to the null geometry of the
problem. Besides their role in the present problem, we expect these
estimates to be useful for other derivative-type wave equations, including
coupled systems with multiple propagation speeds.
These weighted energy and local energy estimates are stated in
\cref{thm:weighted-energy} below.

The rest of the paper is organized as follows.
In Section~\ref{sec:preparation}, we introduce the null-coordinate notation,
collect the one-dimensional inequalities used throughout the paper, and recall
the local existence result.
In Section~\ref{sec:weighted-energy}, we establish
\cref{thm:weighted-energy} and prove the corresponding radial pointwise
decay estimate in \cref{lem:radial-pointwise-estimate}.
In Section~\ref{sec:proof-main-theorem}, we combine these estimates with a
bootstrap argument to prove \cref{thm:main}.

\section{Preparation}\label{sec:preparation}
\subsection{Notation}

The proof of the main theorem is based primarily on multiplier arguments in
spacetime null coordinates. We begin by introducing the notation and recording
several basic facts needed later.

As usual, the notation $a\lesssim b$ (or equivalently $b\gtrsim a$)
means that there exists a constant $C>0$ such that $a\leq Cb$,
where $C$ may vary from line to line.
Similarly, $a\approx b$ indicates that both $a\lesssim b$
and $b\lesssim a$ hold.

We work mainly in the $(u,v,\Theta)$ coordinates, where
\begin{Eq*}
    u=&\frac{t-r}{2},\quad v=\frac{t+r}{2},\quad  r=|x|,\qquad\kl(\Leftrightarrow~  t=u+v,\quad r=v-u\kr),
\end{Eq*}
and $\Theta=\{\theta^A\}_{A=1}^{n-1}$ denotes the standard spherical coordinates
on $\rS^{n-1}$, so that $v\geq |u|$ for $(t,x)\in \rR_+\times\rR^n$.
In these coordinates, the Minkowski metric takes the form
\begin{Eq*}
    {\bf m}=&-\d t^2+\d x^2=-4\d u\d v+r^2\d\Theta^2,\\
     \d\Theta^2:=&\sm_{AB}\d\theta^A\d\theta^B=(\d\theta^1)^2+\sum_{A=2}^{n-1}\kl(\prod_{D=1}^{A-1}\sin^2\theta^D\kr)(\d\theta^A)^2.
\end{Eq*}
For a scalar field $f$, we write
\begin{Eq*}
    \spartial f:=r^{-1}\nabla_{\rS^{n-1}}f,
    \qquad
    |\spartial f|^2:=r^{-2}\sm^{AB}\partial_Af\partial_Bf
\end{Eq*}
for its angular derivative and the corresponding squared norm.
Throughout the paper, we use the Einstein summation convention and raise and
lower indices with respect to the metric ${\bf m}$. In particular, for any
vector field $X$, we have $X_u=-2X^v$ and $X_v=-2X^u$. The nonzero
Christoffel symbols in the $(u,v,\Theta)$ coordinates are
\begin{Eq*}
    \Gamma^{A}{}_{vD}=-\Gamma^A{}_{uD}=r^{-1}\delta_D^A,\quad \Gamma^v{}_{BD}=-\Gamma^u{}_{BD}=-\frac{1}{2}r^{-1}m_{BD},\quad \Gamma^{A}{}_{BD}=\sGamma^A{}_{BD},
\end{Eq*}
where $\sGamma$ denotes the Christoffel symbols of the spherical metric $\sm$
in the $\Theta$ coordinates. We also introduce the following notation for null
hypersurfaces and spacetime regions:
\begin{Eq*}
S_0^T:=&[0,T]\times\rR^n,\quad&
\rS_{U}^V:=&\{(t,x):u=U,~v=V\},\\
\rH_{U}^{V_1,V_2}:=&\{(t,x):u=U,~V_1\leq v\leq V_2\},\quad&
\rHb_{U_1,U_2}^{V}:=&\{(t,x):U_1\leq u\leq U_2,~v=V\},\\
\rD_{U_1,U_2}^{V_1,V_2}:=&\{(t,x):U_1\leq u\leq U_2,~V_1\leq v\leq V_2,~t\geq 0\}.
\end{Eq*}
The following figure illustrates these surfaces and regions.
\begin{figure}[H]
\centering
\includegraphics[width=0.3\textwidth]{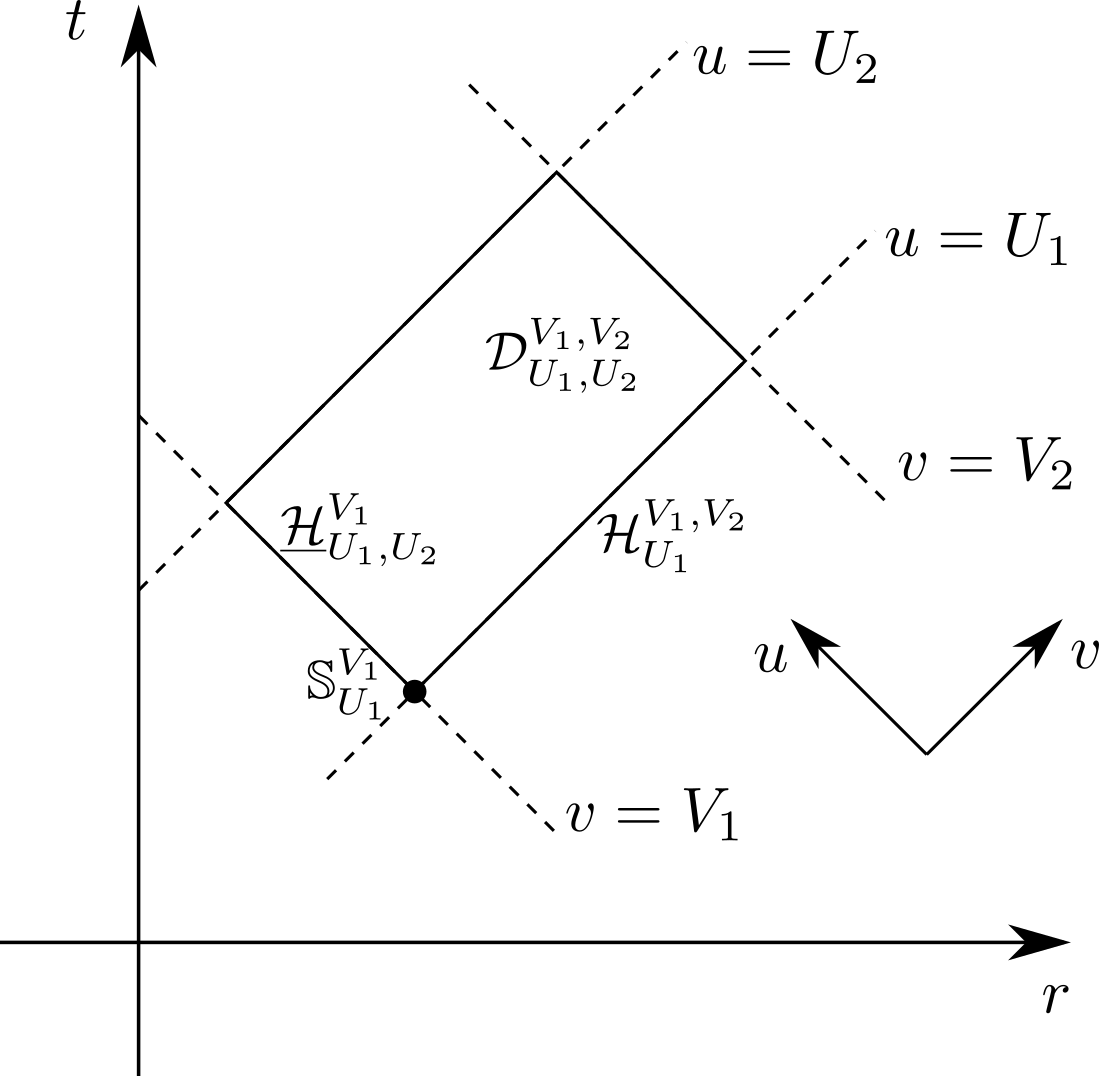}
\end{figure}

We shall also use the scaling vector field and the Lorentz boosts. More
precisely, let
\begin{Eq}\label{eq:commutation-vector-fields}
S:=t\partial_t+r\partial_r,\qquad
L_i:=t\partial_i+x_i\partial_t,\quad 1\leq i\leq n,\qquad
Z:=\{\partial_t,\partial_{x_1},\ldots,\partial_{x_n},S,L_1,\ldots,L_n\}.
\end{Eq}
For a radial scalar field, we also write
\begin{Eq*}
L:=t\partial_r+r\partial_t.
\end{Eq*}

\subsection{One-dimensional inequalities and local existence}

We record the standard one-dimensional
inequalities; see, e.g., \cite[Lemma 6.1]{DaiMeiWeiYang2025}.
\begin{lemma}\label{lem:one-dimensional-weighted-estimates}
Let $h\in C^1([s_1,s_2])$ and
$g\in C([s_1,s_2])\cap C^1((s_1,s_2])$, where $g$ is nonnegative and
increasing. Assume that the quantities on the right-hand sides below are
finite. Then
\begin{Eq}\label{eq:one-dimensional-weighted-first}
|h(s_2)|^2g(s_2)^{-1}
+g(s_2)^{-2}\int_{s_1}^{s_2}|h(s)|^2g'(s)\d s
\lesssim
|h(s_1)|^2g(s_2)^{-1}
+\int_{s_1}^{s_2}|h'(s)|^2g'(s)^{-1}\d s,
\end{Eq}
\begin{Eq}\label{eq:one-dimensional-weighted-second}
|h(s_2)|^2g(s_2)^{-1}
+\int_{s_1}^{s_2}|h(s)|^2g(s)^{-2}g'(s)\d s
\lesssim
|h(s_1)|^2g(s_1)^{-1}
+\int_{s_1}^{s_2}|h'(s)|^2g'(s)^{-1}\d s.
\end{Eq}
\end{lemma}

We use the local well-posedness theory for radial $H^2\times H^1$ data.
An adaptation of the contraction arguments in
\cite[Sections~4--5]{MR2980460} to the present system gives the following result.
\begin{proposition}\label{prop:local-well-posedness}
Let $n\geq3$, $1<p,q<n/(n-2)$, and
\begin{Eq*}
(f_\phi,g_\phi),\ (f_\psi,g_\psi)
\in H^2_{\mathrm{rad}}(\rR^n)\times H^1_{\mathrm{rad}}(\rR^n).
\end{Eq*}
Here \(H_{\mathrm{rad}}^s(\mathbb R^n)\) denotes the subspace of \(H^s(\mathbb R^n)\) consisting of radial functions.
For every $\varepsilon>0$, there exists $T_{\mathrm{loc}}>0$ such that
\eqref{eq:main-system} admits a unique solution satisfying
\begin{Eq*}
\phi,\psi\in C\kl([0,T_{\mathrm{loc}}];H^2_{\mathrm{rad}}(\rR^n)\kr)
\cap C^1\kl([0,T_{\mathrm{loc}}];H^1_{\mathrm{rad}}(\rR^n)\kr).
\end{Eq*}
The time $T_{\mathrm{loc}}$ can be chosen uniformly for Cauchy data in
bounded subsets of this Sobolev space. Solutions satisfy the
domain-of-dependence property and, in particular, have finite propagation speed.
Moreover, a solution on $[0,T)$, with $T<\infty$, extends beyond $T$ whenever
\begin{Eq*}
\sup_{0\leq t<T}\kl(
\|\phi(t)\|_{H^2}+\|\phi_t(t)\|_{H^1}
+\|\psi(t)\|_{H^2}+\|\psi_t(t)\|_{H^1}
\kr)<\infty.
\end{Eq*}
\end{proposition}

\section{Weighted energy and pointwise estimates}\label{sec:weighted-energy}

In this section, we establish the weighted energy and local energy estimates
used in the proof of \cref{thm:main}, together with the corresponding
pointwise decay estimates for radial scalar fields.

We first state the results and then present their proofs in the following
subsections.
To begin with, we set $c_n:=(n-1)(n-3)/4$ and introduce the modified
derivatives as follows. For any vector field $X$, we define
\begin{Eq*}
\hat Xf:=Xf+\frac{n-1}{2}\frac{Xr}{r} f
=r^{-(n-1)/2}X\kl(r^{(n-1)/2}f\kr).
\end{Eq*}
In particular,
\begin{Eq*}
\hpartial{v}f:=\partial_vf+\frac{n-1}{2r}f,
\qquad
\hpartial{u}f:=\partial_uf-\frac{n-1}{2r}f,
\qquad
\hpartial{r}f:=\partial_rf+\frac{n-1}{2r}f.
\end{Eq*}
With this notation, we have the following weighted energy and local energy
estimate.
We omit the natural volume elements from spacetime and hypersurface
integrals unless ambiguity may arise.
For example, for $f=f(u,v)$, the following identities hold whenever
the integration ranges lie entirely in $v\geq|u|$:
\begin{Eq*}
    \int_{\rH_{u}^{v_1,v_2}}f=|\rS^{n-1}|\int_{v_1}^{v_2}f(u,v)r^{n-1}\d v,\quad \int_{\rD_{u_1,u_2}^{v_1,v_2}}f=2|\rS^{n-1}|\int_{v_1}^{v_2}\int_{u_1}^{u_2}f(u,v)r^{n-1}\d u\d v.
\end{Eq*}

\begin{theorem}[Weighted energy and local energy estimates]
\label{thm:weighted-energy}
Let $n\geq3$. The estimate
\begin{Eq}\label{eq:E2-estimate}
E_2[f](T)
\lesssim E_0[f]
+\int_{S_0^T}|\square f|
\kl(\kl<u\kr>|\hpartial{u}f|+\kl<v\kr>|\hpartial{v}f|\kr)
\end{Eq}
holds for any $T>0$ and any smooth real-valued scalar field $f$ on $S_0^T$
with sufficient spatial decay, where 
\begin{Eq}\label{eq:energy-definitions}
E_2[f](T)
:={}&E_1[f](T)
+\sup_{j\in\rZ}
\int_{S_0^T\cap\{2^{j-1}\leq r\leq2^j\}}r^{-1}Q[f],\\
E_1[f](T)
:={}&\sup_{u\in(-\infty,T/2]}
\int_{\rH_u^{|u|,T-u}}
\kl<v\kr>|\hpartial{v}f|^2
+\kl<u\kr>\kl(|\spartial f|^2+c_n|r^{-1}f|^2\kr)\\
&+\sup_{v\in[0,\infty)}
\int_{\rHb_{-v,\min\{v,T-v\}}^v}
\kl<u\kr>|\hpartial{u}f|^2
+\kl<v\kr>\kl(|\spartial f|^2+c_n|r^{-1}f|^2\kr)\\
&+\sup_{t\in[0,T]}\int_{\rR^n}Q[f],\\
E_0[f]
:={}&\int_{t=0}\kl<r\kr>
\kl(|\partial_tf|^2+|\nabla_xf|^2\kr),
\end{Eq}
and
\begin{Eq}\label{eq:Q-definition}
Q[f]:={}&\kl<u\kr>|\hpartial{u}f|^2
+\kl<v\kr>|\hpartial{v}f|^2
+\kl<v\kr>\kl(|\spartial f|^2+c_n|r^{-1}f|^2\kr).
\end{Eq}
The implicit constant in \eqref{eq:E2-estimate} depends only on $n$.
\end{theorem}

As a consequence, we obtain the following weighted spacetime estimate,
analogous to the classical Keel--Smith--Sogge estimates \cite{MR1945285}.
\begin{corollary}\label{cor:weighted-spacetime}
Under the assumptions of \cref{thm:weighted-energy}, for every $\delta>0$
and $\eta\leq1$, we have
\begin{Eq}\label{eq:weighted-spacetime-source-estimate}
A_\eta(T)^{-1}
\int_{S_0^T}r^{-1+\delta}\kl<r\kr>^{-\delta+\eta}Q[f]
\lesssim{}&E_0[f]+\int_{S_0^T}|\square f|
\kl(\kl<u\kr>|\hpartial{u}f|+\kl<v\kr>|\hpartial{v}f|\kr),
\end{Eq}
where
\begin{Eq}\label{eq:A-eta-definition}
A_\eta(T):=
\begin{cases}
1,&\eta<0,\\
\ln(2+T),&\eta=0,\\
\kl<T\kr>^\eta,&\eta>0.
\end{cases}
\end{Eq}
The implicit constant depends only on $n,\delta,\eta$.
\end{corollary}

We also establish the following pointwise decay estimate of
Klainerman--Sobolev type for radial scalar fields, expressed in terms
of the weighted energies above; see \cite{MR865359} for the classical
inequalities. It will be used to control the nonlinear terms in
Section~\ref{sec:proof-main-theorem}.

Here and below, $Z$ is as in \eqref{eq:commutation-vector-fields}, and we use the convention
$E_i[Z^{\leq 1}f]:=\sum_{|\alpha|\leq 1}E_i[Z^\alpha f]$ for $i=0,1,2$.

\begin{lemma}\label{lem:radial-pointwise-estimate}
Let $n\geq3$, $T>0$, and let $f=f(u,v)$ be a smooth real-valued radial
scalar field on $S_0^T$ with sufficient spatial decay. Then, for
$0<t\leq T$ and $r>0$, we have
\begin{Eq}\label{eq:radial-pointwise-estimate}
|f_r(u,v)|+|f_t(u,v)|&\lesssim
r^{-\frac{n-2}{2}}\kl<u\kr>^{-1}\kl<v\kr>^{-\frac12}
\kl(E_0[Z^{\leq 1}f]+E_1[Z^{\leq 1}f](u+v)\kr)^{1/2}.
\end{Eq}
The implicit constant depends only on $n$.
\end{lemma}

The proofs start with a modified current identity. We then use $X_1$ to
establish the weighted energy estimate and a second family $X_2$ to obtain
the local energy estimate. Next, we prove \cref{cor:weighted-spacetime}
by dyadic summation. Finally, we establish \cref{lem:radial-pointwise-estimate}
using the weighted null fluxes and one-dimensional inequalities.

\subsection{The modified current identity}

Let $f$ be a smooth real-valued scalar field. We define its stress-energy
tensor by
\begin{Eq*}
T^{\alpha\beta}[f]
:=\nabla^\alpha f\nabla^\beta f
-\frac{1}{2}m^{\alpha\beta}\nabla_\gamma f\nabla^\gamma f.
\end{Eq*}
Since $\square=\partial_t^2-\Delta=-\nabla^\alpha\nabla_\alpha$, it satisfies
\begin{Eq*}
\nabla_\alpha T^{\alpha\beta}[f]
=-\square f\cdot \nabla^\beta f.
\end{Eq*}
In the null frame introduced above, the relevant components are
\begin{Eq*}
T^{uu}=&\frac{1}{4}|\partial_vf|^2,
\quad
T^{uv}=\frac{1}{4}|\spartial f|^2,
\quad
T^{vv}=\frac{1}{4}|\partial_uf|^2,
\quad
m_{AB}T^{AB}=|\spartial f|^2
-\frac{n-1}{2}\nabla_\gamma f\nabla^\gamma f.
\end{Eq*}

For vector fields $X,Y$ and a scalar function $\chi$, we define the modified
current by
\begin{Eq*}
J^\alpha[f;X,\chi,Y]
:=-T^{\alpha\beta}[f]X_\beta
+\frac{1}{2}f^2\nabla^\alpha\chi
-\chi f\nabla^\alpha f
+Y^\alpha.
\end{Eq*}
A direct computation gives
\begin{Eq}\label{eq:modified-current-divergence}
\nabla_\alpha J^\alpha
=&\square f\,(Xf+\chi f)
-T^{\alpha\beta}\pi_{\alpha\beta}
-\chi\nabla_\alpha f\nabla^\alpha f
+\frac{1}{2}f^2\nabla_\alpha\nabla^\alpha\chi
+\nabla_\alpha Y^\alpha,
\end{Eq}
where $\pi$ is the deformation tensor of $X$, defined by
\begin{Eq*}
\pi_{\alpha\beta}
:=&\frac{1}{2}\kl(\nabla_\alpha X_\beta+\nabla_\beta X_\alpha\kr)\\
=&\frac{1}{2}\Big(
m_{\beta\gamma}\kl(\partial_\alpha X^\gamma
+\Gamma^\gamma{}_{\alpha\sigma}X^\sigma\kr)
+m_{\alpha\gamma}\kl(\partial_\beta X^\gamma
+\Gamma^\gamma{}_{\beta\sigma}X^\sigma\kr)
\Big).
\end{Eq*}
Integrating over a spacetime region $\rD$ and
applying the divergence theorem yields
the energy identity
\begin{Eq}\label{eq:divergence-theorem}
\int_{\partial\rD}J^\alpha\bn_\alpha
=\int_{\rD}\nabla_\alpha J^\alpha,
\end{Eq}
where $\bn$ denotes the outward conormal to $\partial\rD$.

We now consider vector fields of the form
\begin{Eq*}
X=X^u\partial_u+X^v\partial_v,
\end{Eq*}
where $X^u$ and $X^v$ are functions of $u$ and $v$. The nonzero
components of the deformation tensor are
\begin{Eq*}
\pi_{uu}=-2\partial_uX^v,
\quad
\pi_{uv}=-\partial_vX^v-\partial_uX^u,
\quad\pi_{vv}=-2\partial_vX^u,
\quad
\pi_{AB}=r^{-1}(X^v-X^u)m_{AB}.
\end{Eq*}
We then choose
\begin{Eq*}
\chi:={}&\frac{n-1}{2r}(X^v-X^u),\\
Y:={}&\frac{n-1}{4r}f^2
\kl(\partial_uX^v\partial_v-\partial_vX^u\partial_u\kr)\\
&+\frac{n-1}{8r^{n-1}}
\kl(
\partial_v\kl(r^{n-2}(X^u+X^v)f^2\kr)\partial_u
-\partial_u\kl(r^{n-2}(X^u+X^v)f^2\kr)\partial_v
\kr).
\end{Eq*}
With these choices, the null components of the current are
\begin{Eq*}
J^u
=&-T^{uu}X_u-T^{uv}X_v
+\frac{1}{2}f^2\nabla^u\chi
-\chi f\nabla^u f+Y^u\\
=&\frac{X^v}{2}|\partial_v f|^2
+\frac{X^u}{2}|\spartial f|^2
-\frac{1}{4}f^2\partial_v\kl(\frac{n-1}{2r}(X^v-X^u)\kr)
+\frac{n-1}{4r}(X^v-X^u)f\partial_vf
-\frac{n-1}{4r}f^2\partial_vX^u\\
&+\frac{n-1}{8r^{n-1}}
\partial_v\kl(r^{n-2}(X^u+X^v)f^2\kr)\\
=&\frac{X^v}{2}|\hpartial{v}f|^2
+\frac{X^u}{2}\kl(|\spartial f|^2+c_n|r^{-1}f|^2\kr).
\end{Eq*}
By symmetry, the other null component is
\begin{Eq*}
J^v
=&\frac{X^u}{2}|\hpartial{u}f|^2
+\frac{X^v}{2}\kl(|\spartial f|^2+c_n|r^{-1}f|^2\kr).
\end{Eq*}
We also have $J^t=J^u+J^v$.

To compute the divergence term, we first suppose that $X^v=0$, so that
$X=X^u\partial_u$. Substitution into
\eqref{eq:modified-current-divergence} gives
\begin{Eq*}
\nabla_\alpha J^\alpha
=&-T^{\alpha\beta}\pi_{\alpha\beta}
-\chi\nabla_\alpha f\nabla^\alpha f
+\frac{1}{2}f^2\nabla_\alpha\nabla^\alpha\chi
+\nabla_\alpha Y^\alpha
+\square f\,(Xf+\chi f)\\
=&\frac{\partial_vX^u}{2}|\partial_uf|^2
+\frac{\partial_uX^u}{2}|\spartial f|^2
+r^{-1}X^u\kl(
|\spartial f|^2
-\frac{n-1}{2}\nabla_\alpha f\nabla^\alpha f
\kr)\\
&\quad
+\frac{n-1}{2r}X^u\nabla_\alpha f\nabla^\alpha f
+\frac{1}{2}f^2\partial_u\partial_v
\kl(\frac{n-1}{2r}X^u\kr)+\frac{n-1}{4r}f^2(\partial_u-\partial_v)
\kl(\frac{n-1}{2r}X^u\kr)\\
&\quad
-\partial_u\kl(\frac{n-1}{4r}f^2\partial_vX^u\kr)
+\frac{n-1}{r}\kl(\frac{n-1}{4r}f^2\partial_vX^u\kr)\\
&\quad
+\frac{n-1}{8r^{n-1}}
(\partial_u\partial_v-\partial_v\partial_u)
\kl(r^{n-2}X^uf^2\kr)
+\square f\kl(X^u\partial_uf-\frac{n-1}{2r}X^uf\kr)\\
=&\frac{\partial_vX^u}{2}|\hpartial{u}f|^2
+\kl(\frac{\partial_uX^u}{2}+\frac{X^u}{r}\kr)
\kl(|\spartial f|^2+c_n|r^{-1}f|^2\kr)
+\square f\,(X^u\hpartial{u}f).
\end{Eq*}
By symmetry, for a general $X=X^u(u,v)\partial_u+X^v(u,v)\partial_v$ we obtain
\begin{Eq}\label{eq:modified-current-identity}
\nabla_\alpha J^\alpha
=&\kl(\nabla_\alpha J^\alpha\kr)_{\mathrm{bulk}}
+\square f\,(X^v\hpartial{v}f+X^u\hpartial{u}f),\\
\kl(\nabla_\alpha J^\alpha\kr)_{\mathrm{bulk}}
:=&\frac{\partial_uX^v}{2}|\hpartial{v}f|^2
+\frac{\partial_vX^u}{2}|\hpartial{u}f|^2+\kl(
\frac{\partial_vX^v+\partial_uX^u}{2}
-\frac{X^v-X^u}{r}
\kr)
\kl(|\spartial f|^2+c_n|r^{-1}f|^2\kr).
\end{Eq}

\subsection{The weighted energy estimate}

For the weighted energy estimate, we choose
\begin{Eq*}
X_1^u:=1+|u|,
\qquad
X_1^v:=1+v.
\end{Eq*}
Here we note that, to the best of our knowledge, multipliers of the form
$|u|^a\partial_u+v^a\partial_v$, $1\leq a<2$, first appeared in the work
of Lindblad and Sterbenz \cite{MR2253534}, who established the associated
deformation-tensor positivity in this range.
For our choice of $X_1$, the null components of $J_1$ and the bulk part of
its divergence are
\begin{Eq*}
J_1^u=&\frac{1+v}{2}|\hpartial{v}f|^2
+\frac{1+|u|}{2}
\kl(|\spartial f|^2+c_n|r^{-1}f|^2\kr),\\
J_1^v=&\frac{1+|u|}{2}|\hpartial{u}f|^2
+\frac{1+v}{2}
\kl(|\spartial f|^2+c_n|r^{-1}f|^2\kr),\\
\kl(\nabla_\alpha J_1^\alpha\kr)_{\mathrm{bulk}}
=&-{\bf 1}_{\{u<0\}}\frac{t}{r}
\kl(|\spartial f|^2+c_n|r^{-1}f|^2\kr),
\end{Eq*}
where ${\bf 1}_A$ denotes the indicator function of $A$.

\begin{proof}[Proof of the energy estimate in \cref{thm:weighted-energy}]
Notice that $J_1^u$ and $J_1^v$ are nonnegative. Applying
\eqref{eq:divergence-theorem} to $J_1$ on $S_0^T$ and its null-truncated
subregions,
discarding the nonpositive term
$\kl(\nabla_\alpha J_1^\alpha\kr)_{\mathrm{bulk}}$, we obtain
\begin{Eq*}
&\sup_{u\in(-\infty,T/2]}
\int_{\rH_u^{|u|,T-u}}J_1^u+\sup_{v\in[0,\infty)}
\int_{\rHb_{-v,\min\{v,T-v\}}^v}J_1^v
+\sup_{t\in[0,T]}\int_{\rR^n}J_1^t\\
\lesssim&\int_{t=0}J_1^t+\int_{S_0^T}|\square f|
\kl|(1+|u|)\hpartial{u}f+(1+v)\hpartial{v}f\kr|.
\end{Eq*}
Since $v\geq|u|$ on $S_0^T$, the three terms on the left are comparable,
respectively, to the three terms defining $E_1[f](T)$ in
\eqref{eq:energy-definitions}.

We next estimate the first term on the right-hand side. On $t=0$, we have
$u=-r/2$ and $v=r/2$, so
\begin{Eq*}
J_1^t
={}&\kl(1+\frac r2\kr)\kl(
|\partial_tf|^2
+\kl|\partial_rf+\frac{n-1}{2r}f\kr|^2
+|\spartial f|^2+c_n r^{-2}|f|^2\kr).
\end{Eq*}
Meanwhile, the classical Hardy inequality and integration by parts on $t=0$ give
\begin{Eq*}
\left\|\kl<r\kr>^{\frac12}r^{-1}f\right\|_{L^2}^2
\lesssim&\left\|\nabla_x\kl(\kl<r\kr>^{\frac12}f\kr)\right\|_{L^2}^2\\
={}&\int_{\rR^n}\kl<r\kr>|\nabla_xf|^2\d x-\int_{\rR^n}\kl<r\kr>^{\frac12}
\Delta_x\kl(\kl<r\kr>^{\frac12}\kr)|f|^2\d x\\
\leq{}&\int_{\rR^n}\kl<r\kr>|\nabla_xf|^2\d x,
\end{Eq*}
since $\Delta_x\kl(\kl<r\kr>^{\frac12}\kr)\geq0$.
Consequently,
\begin{Eq*}
\int_{t=0}J_1^t\lesssim E_0[f].    
\end{Eq*}

Substituting these comparisons into the flux estimate, we conclude that
\begin{Eq}\label{eq:E1-estimate}
E_1[f](T)
\lesssim E_0[f]
+\int_{S_0^T}
|\square f|\,
\kl|
(1+|u|)\hpartial{u}f
+(1+v)\hpartial{v}f
\kr|.
\end{Eq}
This proves the energy part of \cref{thm:weighted-energy}.
Here we remark that, since $X_1$ is only piecewise $C^1$, we in fact apply
the divergence theorem separately in the regions $u\geq0$ and $u\leq0$,
while the fluxes across $u=0$ cancel by continuity of $J_1^u$.
The same argument also applies to the second multiplier below.
\end{proof}

\subsection{The weighted local energy estimate}
\label{subsec:weighted-energy-estimates}

To derive the weighted local energy estimate, for any fixed $R>0$, we choose
\begin{Eq*}
h_R:=\kl(\frac{r}{r+R}\kr)^{\frac12},
\qquad
X_2^u:=-(1+|u|)h_R,
\qquad
X_2^v:=
\begin{cases}
(1+3v)h_R,&u\geq0,\\
(1+3r)h_R,&u<0.
\end{cases}
\end{Eq*}
A radial weight of this form was also used in Hidano--Wang--Yokoyama
\cite[(2.14)]{HidanoWangYokoyama2012}. For further background, see that
paper and the references therein.
For this choice, the null components of $J_2$ in the region $u\geq0$ are
\begin{Eq*}
J_2^u
={}&\frac{1+3v}{2}h_R|\hpartial{v}f|^2
-\frac{1+u}{2}h_R
\kl(|\spartial f|^2+c_n|r^{-1}f|^2\kr),\\
J_2^v
={}&-\frac{1+u}{2}h_R|\hpartial{u}f|^2
+\frac{1+3v}{2}h_R
\kl(|\spartial f|^2+c_n|r^{-1}f|^2\kr).
\end{Eq*}
Since $h_R\leq1$, $J_2^u$ and $J_2^v$ are controlled by
$J_1^u$ and $J_1^v$ in the region $u\geq 0$, respectively.
Moreover, identity \eqref{eq:modified-current-identity} yields
\begin{Eq*}
\kl(\nabla_\alpha J_2^\alpha\kr)_{\mathrm{bulk}}
={}&-\frac{1+3v}{2}h_R'|\hpartial{v}f|^2
-\frac{1+u}{2}h_R'|\hpartial{u}f|^2\\
&-\kl(
(2+u+3v)\kl(\frac{h_R}{r}-\frac{h_R'}{2}\kr)-h_R
\kr)
\kl(|\spartial f|^2+c_n|r^{-1}f|^2\kr).
\end{Eq*}
Here $h_R'\geq 0$ and
$h_R'\leq h_R/(2r)$, so
\begin{Eq*}
\frac{1+3v}{2}h_R'
\gtrsim\kl<v\kr>h_R'\geq0,&\qquad
\frac{1+u}{2}h_R'
\gtrsim\kl<u\kr>h_R'\geq0,\\
(2+u+3v)\kl(\frac{h_R}{r}-\frac{h_R'}{2}\kr)-h_R
&\geq\frac{6+7u+5v}{4r}h_R
\gtrsim\frac{\kl<v\kr>}{r}h_R\geq0.
\end{Eq*}

In the region $u<0$, we similarly obtain
\begin{Eq*}
J_2^u
={}&\frac{1+3r}{2}h_R|\hpartial{v}f|^2
-\frac{1-u}{2}h_R
\kl(|\spartial f|^2+c_n|r^{-1}f|^2\kr),\\
J_2^v
={}&-\frac{1-u}{2}h_R|\hpartial{u}f|^2
+\frac{1+3r}{2}h_R
\kl(|\spartial f|^2+c_n|r^{-1}f|^2\kr).
\end{Eq*}
Again, since $h_R\leq1$ and $1+3r\lesssim1+v$ in this region,
$J_2^u$ and $J_2^v$ are controlled by
$J_1^u$ and $J_1^v$ in the region $u< 0$, respectively.
Again,  \eqref{eq:modified-current-identity}  shows
\begin{Eq*}
\kl(\nabla_\alpha J_2^\alpha\kr)_{\mathrm{bulk}}
={}&-\kl(\frac32h_R+\frac{1+3r}{2}h_R'\kr)
|\hpartial{v}f|^2
-\frac{1-u}{2}h_R'|\hpartial{u}f|^2\\
&-\kl(
(2+3r-u)\kl(\frac{h_R}{r}-\frac{h_R'}{2}\kr)-2h_R
\kr)
\kl(|\spartial f|^2+c_n|r^{-1}f|^2\kr).
\end{Eq*}
Here, similarly,
\begin{Eq*}
\frac32h_R+\frac{1+3r}{2}h_R'
\gtrsim\kl<v\kr>h_R'\geq0,&\qquad
\frac{1-u}{2}h_R'
\gtrsim\kl<u\kr>h_R'\geq0,\\
(2+3r-u)\kl(\frac{h_R}{r}-\frac{h_R'}{2}\kr)-2h_R
&\geq\frac{6+r-3u}{4r}h_R
\gtrsim\frac{\kl<v\kr>}{r}h_R\geq0.
\end{Eq*}
Thus, since $c_n\geq0$ for $n\geq3$,
all terms in $\kl(\nabla_\alpha J_2^\alpha\kr)_{\mathrm{bulk}}$ are nonpositive in both regions.
Compared with the choice $h_R=r/(r+R)$, the square-root weight has the
advantage of ensuring the required sign of all bulk terms directly.
Alternatively, the same local energy estimate can be obtained by combining
the multiplier identities corresponding to $h_R=1$ and $h_R=r/(r+R)$.

\begin{proof}[Proof of the local energy estimate in \cref{thm:weighted-energy}]
Applying \eqref{eq:divergence-theorem} to $J_2$ on $S_0^T$, we obtain
\begin{Eq*}
-\int_{S_0^T}\kl(\nabla_\alpha J_2^\alpha\kr)_{\mathrm{bulk}}
\leq{}&\int_{t=0}\kl|J_2^t\kr|
+\int_{t=T}\kl|J_2^t\kr|+\int_{S_0^T}|\square f|
\kl|X_2^u\hpartial{u}f+X_2^v\hpartial{v}f\kr|.
\end{Eq*}
We first estimate the left-hand side. On the annulus $R/2<r<2R$, we have
$h_R\approx1$ and
$h_R'\approx r^{-1}$.
Hence, on this annulus, the coefficient bounds above yield
\begin{Eq*}
-\kl(\nabla_\alpha J_2^\alpha\kr)_{\mathrm{bulk}}
\gtrsim{}&\frac{\kl<u\kr>}{r}|\hpartial{u}f|^2
+\frac{\kl<v\kr>}{r}|\hpartial{v}f|^2+\frac{\kl<v\kr>}{r}
\kl(|\spartial f|^2+c_n|r^{-1}f|^2\kr),
\end{Eq*}
in both regions $u\geq0$ and $u<0$.
It therefore follows that
\begin{Eq*}
\int_{S_0^T\cap\{R/2<r<2R\}}r^{-1}Q[f]
\lesssim-\int_{S_0^T}\kl(\nabla_\alpha J_2^\alpha\kr)_{\mathrm{bulk}},
\end{Eq*}
with $Q[f]$ defined in \eqref{eq:Q-definition}.

We next estimate the right-hand side. The bounds for $J_2^u$ and $J_2^v$
above give
\begin{Eq*}
\int_{t=0}\kl|J_2^t\kr|
+\int_{t=T}\kl|J_2^t\kr|
\lesssim\int_{t=0}J_1^t
+\int_{t=T}J_1^t
\lesssim E_1[f](T).
\end{Eq*}
Moreover, $h_R\leq1$ and $r\leq2v$ when $u<0$, so the source term satisfies
\begin{Eq*}
|\square f|\kl|X_2^u\hpartial{u}f+X_2^v\hpartial{v}f\kr|
\lesssim|\square f|
\kl(\kl<u\kr>|\hpartial{u}f|+\kl<v\kr>|\hpartial{v}f|\kr).
\end{Eq*}
Finally, setting $R=2^j$, taking the supremum over $j\in\rZ$, and using
\eqref{eq:E1-estimate}, we conclude that
\begin{Eq*}
E_2[f](T)
\lesssim E_0[f]
+\int_{S_0^T}|\square f|
\kl(\kl<u\kr>|\hpartial{u}f|+\kl<v\kr>|\hpartial{v}f|\kr).
\end{Eq*}
This completes the proof of \cref{thm:weighted-energy}.
\end{proof}

\subsection{Proof of the KSS-type estimate}

\begin{proof}[Proof of \cref{cor:weighted-spacetime}]
We first establish
\begin{Eq}\label{eq:E2-weighted-spacetime}
A_\eta(T)^{-1}
\int_{S_0^T}r^{-1+\delta}\kl<r\kr>^{-\delta+\eta}Q[f]
\lesssim E_2[f](T).
\end{Eq}
This follows by a standard dyadic summation argument; see, for example,
\cite[the proof of (3.5)]{MR2980460}. We first consider $r\leq1$. Since
$\kl<r\kr>\approx1$ and $\delta>0$,
\begin{Eq*}
\int_{S_0^T\cap\{r\leq1\}}
r^{-1+\delta}\kl<r\kr>^{-\delta+\eta}Q[f]
\lesssim&
\sum_{j\leq0}2^{\delta j}
\int_{S_0^T\cap\{2^{j-1}\leq r\leq2^j\}}
r^{-1}Q[f]\\
\lesssim&
\sum_{j\leq0}2^{\delta j}E_2[f](T)\\
\lesssim& E_2[f](T).
\end{Eq*}

For $1<r\leq T+1$, we have $\kl<r\kr>\approx r$, and hence
\begin{Eq*}
\int_{S_0^T\cap\{1<r\leq T+1\}}
r^{-1+\delta}\kl<r\kr>^{-\delta+\eta}Q[f]
\lesssim&
\sum_{1\leq2^{j-1}\leq T+1}2^{\eta j}
\int_{S_0^T\cap\{2^{j-1}\leq r\leq2^j\}}
r^{-1}Q[f]\\
\lesssim&
\sum_{1\leq j\lesssim\ln(T+2)}2^{\eta j}E_2[f](T)\\
\lesssim& A_\eta(T)E_2[f](T).
\end{Eq*}

It remains to consider $r>T+1$. We split this region as
\begin{Eq*}
\{r>T+1\}
=&\{r>T+1,\ r\leq2(t+1)\}\cup\{r>T+1,\ r>2(t+1)\}.
\end{Eq*}
On the first region, we have $r\approx\kl<r\kr>\approx \kl<t\kr>\approx\kl<T\kr>$.
The dyadic term in
$E_2[f](T)$ directly gives
\begin{Eq*}
\int_{S_0^T\cap\{r>T+1,\ r\leq2(t+1)\}}
r^{-1+\delta}\kl<r\kr>^{-\delta+\eta}Q[f]
\lesssim& A_\eta(T)E_2[f](T).
\end{Eq*}

On the second region, $r>2(t+1)$, the integral of $Q[f]$ is controlled by
the time-slice term of $E_1[f](T)$ in \eqref{eq:energy-definitions}.
Since $\eta\leq1$ and
$r>T+1$, we obtain
\begin{Eq*}
&\int_{S_0^T\cap\{r>T+1,\ r>2(t+1)\}}
r^{-1+\delta}\kl<r\kr>^{-\delta+\eta}Q[f]\\
&\qquad\lesssim
\kl<T\kr>^{-1+\eta}
\int_0^T
\int_{\{r>2(t+1)\}}Q[f]\\
&\qquad\lesssim
T\kl<T\kr>^{-1+\eta}E_1[f](T)\\
&\qquad\lesssim A_\eta(T)E_2[f](T).
\end{Eq*}
Combining the three regions and using $A_\eta(T)\gtrsim1$, we obtain
\eqref{eq:E2-weighted-spacetime}. Applying \eqref{eq:E2-estimate} gives
the asserted bound.
\end{proof}

\subsection{Proof of the pointwise decay estimate}

\begin{proof}[Proof of \cref{lem:radial-pointwise-estimate}]
We first establish a weak radial estimate. Fix a point $(U,V)$ with
$V>|U|$ and $U+V\leq T$.

Suppose first that $U\geq0$. Applying
\eqref{eq:one-dimensional-weighted-first} with
$h(v)=(v-U)^{\frac{n-1}{2}}f(U,v)$ and $g(v)=v-U$ on $[U,V]$, we obtain
\begin{Eq*}
(V-U)^{n-2}|f(U,V)|^2
\lesssim
\int_U^V(v-U)^{n-1}|\hpartial{v}f(U,v)|^2\d v.
\end{Eq*}
Here $h(U)=0$ by the definition of $h$ and the regularity of $f$ at
the axis $r=0$, so the boundary term vanishes.

Suppose next that $U<0$. Applying
\eqref{eq:one-dimensional-weighted-second} with the same $h$ and $g$
on $[-U,V]$, we obtain
\begin{Eq*}
(V-U)^{n-2}|f(U,V)|^2
\lesssim{}&(-2U)^{n-2}|f(U,-U)|^2+\int_{-U}^V(v-U)^{n-1}|\hpartial{v}f(U,v)|^2\d v.
\end{Eq*}
In both cases, the integral term is controlled by the outgoing part of
$E_1[f](U+V)$ in \eqref{eq:energy-definitions}, since $v\geq|U|$ and
\begin{Eq*}
\int_{|U|}^V(v-U)^{n-1}|\hpartial{v}f(U,v)|^2\d v
&\lesssim\kl<U\kr>^{-1}
\int_{|U|}^V\kl<v\kr>|\hpartial{v}f(U,v)|^2(v-U)^{n-1}\d v\\
&\lesssim\kl<U\kr>^{-1}E_1[f](U+V).
\end{Eq*}
For the initial term, weighted Cauchy--Schwarz on the initial hypersurface,
written in $(t,r)$ coordinates, gives
\begin{Eq*}
|f(0,r)|^2
=&\left|\int_r^\infty\partial_s f(0,s)\d s\right|^2\\
\leq{}&\left(\int_r^\infty\kl<s\kr>s^{n-1}
|\partial_s f(0,s)|^2\d s\right)
\left(\int_r^\infty\frac{\d s}{\kl<s\kr>s^{n-1}}\right)\\
\lesssim{}&E_0[f]\,r^{2-n}\kl<r\kr>^{-1}.
\end{Eq*}
For $U<0$, the initial endpoint $(U,-U)$ has
$r=-2U$, so this yields
\begin{Eq*}
(-2U)^{n-2}|f(U,-U)|^2
\lesssim\kl<U\kr>^{-1}E_0[f].
\end{Eq*}
Combining the two cases, we obtain
\begin{Eq}\label{eq:weighted-radial-pointwise-intermediate}
|f(u,v)|
\lesssim&
r^{-\frac{n-2}{2}}\kl<u\kr>^{-\frac12}
\kl(E_0[f]+E_1[f](u+v)\kr)^{1/2}.
\end{Eq}

Introduce the family $\tilde Z:=\{\partial_t,\partial_r,S,L\}$.
Since $f$ is radial, we have
\begin{Eq*}
\partial_k f=\omega_k f_r,
\qquad
L_k f=\omega_kLf,
\qquad 1\leq k\leq n,
\end{Eq*}
where $\omega_k:=x_k/r$. Consequently,
\begin{Eq}\label{eq:radial-commuted-energy-control}
E_i[f_r]+E_i[Lf]\lesssim E_i[Z^{\leq 1}f],
\qquad i=0,1,2.
\end{Eq}
Each field in
$\tilde Z^{\leq 1}f$ is radial, so applying
\eqref{eq:weighted-radial-pointwise-intermediate} to these fields gives
\begin{Eq}\label{eq:commuted-radial-pointwise}
|\tilde Z^{\leq 1}f(U,V)|
\lesssim&
(V-U)^{-\frac{n-2}{2}}\kl<U\kr>^{-\frac12}
\kl(E_0[Z^{\leq 1}f]+E_1[Z^{\leq 1}f](U+V)\kr)^{1/2}.
\end{Eq}

When $|U|\geq2$, we also have
\begin{Eq*}
\{f_t,f_r\}=\frac{f_v\pm f_u}{2},\qquad
f_u=\frac{Sf-Lf}{2u},\qquad
f_v=\frac{Sf+Lf}{2v}.
\end{Eq*}
Combining these identities with
\eqref{eq:weighted-radial-pointwise-intermediate}, we further obtain
\begin{Eq}\label{eq:radial-pointwise-large-u}
|f_t|(U,V)+|f_r|(U,V)
\lesssim&
(V-U)^{-\frac{n-2}{2}}\kl<U\kr>^{-\frac32}
\kl(E_0[Z^{\leq 1}f]+E_1[Z^{\leq 1}f](U+V)\kr)^{1/2}.
\end{Eq}

We next improve this estimate in the region
$\kl<U\kr>\lesssim\kl<V\kr>$. More precisely, assume that
$\kl<V\kr>\geq16\kl<U\kr>$, so that
$V-U\approx V\approx\kl<V\kr>$, and consider
\begin{Eq*}
(u,v)\in\kl[U-\frac{\kl<U\kr>}{8}, U\kr]\times \kl[V-\frac{\kl<U\kr>}{8}, V\kr]\equiv I_U\times I_V.
\end{Eq*}
On this region,
$\kl<u\kr>\approx\kl<U\kr>$ and
$r\approx\kl<v\kr>\approx\kl<V\kr>$. Set
$h:=r^{\frac{n-1}{2}}f$. For any $(u,v)\in I_U\times I_V$, the
fundamental theorem of calculus gives
\begin{Eq*}
h_t(U,V)
={}&h_t(u,v)+\int_u^U\partial_u h_t(u',v)\d u'
+\int_v^V\partial_v h_t(U,v')\d v'.
\end{Eq*}
Taking absolute values, averaging over $I_U\times I_V$, and applying
Cauchy--Schwarz, with $|I_U|=|I_V|=\kl<U\kr>/8$, we obtain
\begin{Eq*}
|h_t(U,V)|
\lesssim{}&\kl<U\kr>^{-1}\|h_t\|_{L^2(I_U\times I_V)}
+\kl<U\kr>^{\frac12}\sup_{v\in I_V}\|\partial_u h_t\|_{L^2(I_U)}+\kl<U\kr>^{\frac12}\sup_{u\in I_U}\|\partial_v h_t\|_{L^2(I_V)}.
\end{Eq*}
For the first term, using
$2h_t=h_u+h_v$, we estimate
\begin{Eq*}
\|h_u\|_{L^2(I_U\times I_V)}^2
\lesssim& \kl<U\kr>^{-1}\int_{I_V}\int_{I_U}\kl<u\kr>|h_u|^2 \d u\d v
\lesssim  \kl<U\kr>^{-1}E_1[Z^{\leq 1}f](U+V)\int_{I_V}\d v,\\
\|h_v\|_{L^2(I_U\times I_V)}^2
\lesssim& \kl<V\kr>^{-1}\int_{I_U}\int_{I_V}\kl<v\kr>|h_v|^2 \d v\d u 
\lesssim  \kl<U\kr>^{-1}E_1[Z^{\leq 1}f](U+V)\int_{I_U}\d u.
\end{Eq*}
Since $|I_V|=|I_U|=\kl<U\kr>/8$ and $2h_t=h_u+h_v$, it follows that
\begin{Eq*}
\|h_t\|_{L^2(I_U\times I_V)}
\lesssim \Big(E_1[Z^{\leq 1}f](U+V)\Big)^{1/2}.
\end{Eq*}
For the second and third terms, suppose first that $|U|\leq2$. Then
$\kl<u\kr>\approx\kl<U\kr>\approx1\lesssim\kl<v\kr>$, and hence
\begin{Eq*}
|r^{-\frac{n-1}{2}}\partial_u h_t|
=&|\hpartial{u} f_t|\lesssim\kl<U\kr>^{-\frac{3}{2}}\kl(\kl<u\kr>|\hpartial{u} Z^{\leq 1}f|^2\kr)^{\frac{1}{2}},\\
|r^{-\frac{n-1}{2}}\partial_v h_t|
=&|\hpartial{v} f_t|\lesssim\kl<U\kr>^{-\frac{3}{2}}\kl(\kl<v\kr>|\hpartial{v} Z^{\leq 1}f|^2\kr)^{\frac{1}{2}}.
\end{Eq*}
When $|U|\geq2$, we instead use the identities
\begin{Eq*}
r^{-\frac{n-1}{2}}\partial_u h_t
=&\frac{\hpartial{u}(Sf-Lf-2f)}{4u}+\frac{\hpartial{u}(Sf+Lf)}{4v}-\frac{n-1}{4ru}f,\\
r^{-\frac{n-1}{2}}\partial_v h_t
=&\frac{\hpartial{v}(Sf-Lf)}{4u}+\frac{\hpartial{v}(Sf+Lf-2f)}{4v}+\frac{n-1}{4rv}f.
\end{Eq*}
Since
$|u|\approx\kl<U\kr>\lesssim v\approx r\approx\kl<V\kr>$,
the preceding identities and \eqref{eq:commuted-radial-pointwise} imply
\begin{Eq*}
|r^{-\frac{n-1}{2}}\partial_u h_t|
\lesssim&\kl<U\kr>^{-\frac{3}{2}}\kl(\kl<u\kr>|\hpartial{u} Z^{\leq 1}f|^2\kr)^{\frac{1}{2}}+\kl<U\kr>^{-\frac{3}{2}}\kl<V\kr>^{-\frac{n}{2}}\kl(E_0[Z^{\leq 1}f]+E_1[Z^{\leq 1}f](U+V)\kr)^{1/2},\\
|r^{-\frac{n-1}{2}}\partial_v h_t|
\lesssim&\kl<U\kr>^{-\frac{3}{2}}\kl(\kl<v\kr>|\hpartial{v} Z^{\leq 1}f|^2\kr)^{\frac{1}{2}}+\kl<U\kr>^{-\frac{1}{2}}\kl<V\kr>^{-\frac{n+2}{2}}\kl(E_0[Z^{\leq 1}f]+E_1[Z^{\leq 1}f](U+V)\kr)^{1/2}.
\end{Eq*}
Therefore, using $|I_U|=|I_V|=\kl<U\kr>/8$ and
$\kl<V\kr>\gtrsim\kl<U\kr>$, a direct calculation yields
\begin{Eq*}
\sup_{v\in I_V}
\|\partial_u h_t\|_{L^2(I_U)}
+
\sup_{u\in I_U}
\|\partial_v h_t\|_{L^2(I_V)}
\lesssim
\kl<U\kr>^{-\frac32}
\kl(E_0[Z^{\leq 1}f]+E_1[Z^{\leq 1}f](U+V)\kr)^{\frac12}.
\end{Eq*}
Substituting these estimates into the preceding local Sobolev estimate gives
\begin{Eq*}
|h_t(U,V)|
\lesssim{}&\kl<U\kr>^{-1}\kl(E_0[Z^{\leq 1}f]+E_1[Z^{\leq 1}f](U+V)\kr)^{1/2}.
\end{Eq*}

On the other hand, since
$f_r=(Sf-tf_t)/r$, the preceding estimate and
\eqref{eq:commuted-radial-pointwise} yield, in
the region $\kl<V\kr>\geq16\kl<U\kr>$,
\begin{Eq*}
|f_r|
&\lesssim
\frac{t}{r}|f_t|+\frac{1}{r}|Sf|\\
&\lesssim
\kl(
r^{-\frac{n-1}{2}}\kl<U\kr>^{-1}
+
r^{-\frac{n}{2}}\kl<U\kr>^{-\frac12}
\kr)
\kl(E_0[Z^{\leq 1}f]+E_1[Z^{\leq 1}f](U+V)\kr)^{1/2}.
\end{Eq*}
Combining these estimates with \eqref{eq:radial-pointwise-large-u}, and
using \eqref{eq:commuted-radial-pointwise} in the remaining bounded region,
we obtain \eqref{eq:radial-pointwise-estimate}.
\end{proof}

\section{Proof of the main theorem}\label{sec:proof-main-theorem}

In this section, we combine the weighted energy and pointwise estimates of
Section~\ref{sec:weighted-energy} with a bootstrap argument to prove
\cref{thm:main}.

Throughout this section, the implicit constants in $\lesssim$ are independent
of $\varepsilon$, $T$, and the constant $M$ introduced below.

\subsection{Bootstrap framework}\label{subsec:bootstrap-framework}
By interchanging the two equations if necessary, we may assume that
\begin{Eq*}
1<p\leq q<\frac{n}{n-2}.
\end{Eq*}
Consequently, \eqref{eq:critical-Gamma} becomes
\begin{Eq*}
\Gamma(p,q)=\frac{q+1}{pq-1}-\frac{n-1}{2}.
\end{Eq*}
For later use, we set
\begin{Eq*}
s_p:=1-\frac{n-1}{2}(p-1).
\end{Eq*}
We decompose $\phi$ as
\begin{Eq*}
\phi=\phi^{\mathrm{hom}}+\phi^{\mathrm{inh}},
\end{Eq*}
where $\phi^{\mathrm{hom}}$ solves the homogeneous
wave equation with the same Cauchy data as $\phi$.
Thus,
\begin{Eq*}
\square \phi^{\mathrm{hom}}=0,
\qquad
\kl(\phi^{\mathrm{hom}},\phi_t^{\mathrm{hom}}\kr)\big|_{t=0}
=\varepsilon\kl(f_\phi,g_\phi\kr),
\qquad
\square \phi^{\mathrm{inh}}=|\psi_t|^p,
\qquad
\kl(
\phi^{\mathrm{inh}},\phi_t^{\mathrm{inh}}
\kr)\big|_{t=0}=0.
\end{Eq*}
A similar decomposition also applies to $\psi$.

We first establish the following initial-energy bounds.
\begin{claim}\label{claim:initial-energy-bound}
Under the assumptions of \cref{thm:main}, namely,
\begin{Eq*}
\|\boldsymbol f\|_{\mathcal X}\leq1,
\qquad
0<\varepsilon\leq\varepsilon_0\ll1,
\end{Eq*}
we have
\begin{Eq*}
E_0[Z^{\leq 1}\phi^{\mathrm{hom}}]
+E_0[Z^{\leq 1}\psi^{\mathrm{hom}}]
\lesssim\varepsilon^2,\qquad
E_0[Z^{\leq 1}\phi^{\mathrm{inh}}]
\lesssim\varepsilon^{2p},\qquad 
E_0[Z^{\leq 1}\psi^{\mathrm{inh}}]
\lesssim\varepsilon^{2q}.
\end{Eq*}
\end{claim}
The proof is deferred to the end of this section.

We shall establish weighted a priori estimates for the Sobolev solutions
supplied by \cref{prop:local-well-posedness}. In the bootstrap argument below,
we use
\begin{Eq*}
E_2[Z^{\leq 1}\phi](T)+E_2[Z^{\leq 1}\psi](T)<\infty
\end{Eq*}
on every compact time interval $[0,T]$ of existence, with both energies
continuous in $T$.

By a standard approximation argument, it suffices to consider compactly
supported radial initial data. The local $H^2$ solution supplied by
\cref{prop:local-well-posedness} remains spatially compactly supported by
finite propagation. The local energy estimates used in its construction,
together with \cref{thm:weighted-energy} and standard regularization, ensure
that the weighted energies above are finite and continuous in $T$.
The a priori bounds below are independent of the support radius.

By \cref{claim:initial-energy-bound} and \eqref{eq:E2-estimate}, we may fix
a sufficiently large constant $M$, independently of $\varepsilon$, such that for any $t\geq 0$, 
\begin{Eq*}
E_2[Z^{\leq 1}\phi^{\mathrm{hom}}](t)+E_2[Z^{\leq 1}\psi^{\mathrm{hom}}](t)\leq\varepsilon^2M^2,
\qquad 
E_0[Z^{\leq 1}\phi^{\mathrm{inh}}]\leq\varepsilon^{2p}M^{2p},
\qquad 
E_0[Z^{\leq 1}\psi^{\mathrm{inh}}]\leq\varepsilon^{2q}M^{2q}.
\end{Eq*}
Let $T_*$ be the maximal existence time of this solution in the Sobolev
class of \cref{prop:local-well-posedness}, and let $T_\varepsilon$ be defined
by \eqref{eq:lifespan-lower-bound}. Define the maximal bootstrap time by
\begin{Eq*}
T_b:=\sup\Big\{
0<T<\min\{T_*,T_\varepsilon\}:
E_2[Z^{\leq 1}\psi^{\mathrm{inh}}](t)\leq4\varepsilon^2M^2
\text{ for every }0\leq t\leq T
\Big\}.
\end{Eq*}
Local well-posedness, \cref{claim:initial-energy-bound}, and the weighted
energy continuity noted above imply that $T_b>0$, after decreasing
$\varepsilon_0$ if necessary. In what follows, all estimates are carried out for an arbitrary
$0<T<T_b$ under the bootstrap assumption
\begin{Eq}\label{eq:bootstrap-assumption}
E_2[Z^{\leq 1}\psi^{\mathrm{inh}}](t)\leq4\varepsilon^2 M^2,
\qquad 0\leq t\leq T.
\end{Eq}
We shall improve this estimate to
\begin{Eq*}
E_2[Z^{\leq 1}\psi^{\mathrm{inh}}](t)\leq2\varepsilon^2 M^2,
\end{Eq*}
for every $0\leq t\leq T$. The continuation argument will be completed after
the required estimates have been established.

\subsection{Control of \texorpdfstring{$E_2[Z^{\leq 1}\phi^{\mathrm{inh}}]$}{E2[Z<=1 phi inh]}}
Assume \eqref{eq:bootstrap-assumption}. Applying \eqref{eq:E2-estimate} to
$Z^{\leq 1}\phi^{\mathrm{inh}}$, we obtain
\begin{Eq*}
E_2[Z^{\leq 1}\phi^{\mathrm{inh}}](T)
\lesssim E_0[Z^{\leq 1}\phi^{\mathrm{inh}}]
+\int_{S_0^T}
|\square Z^{\leq 1}\phi^{\mathrm{inh}}|
\kl(
\kl<u\kr>|\hpartial{u} Z^{\leq 1}\phi^{\mathrm{inh}}|
+\kl<v\kr>|\hpartial{v} Z^{\leq 1}\phi^{\mathrm{inh}}|
\kr).
\end{Eq*}
Since
\begin{Eq*}
|\square Z^{\leq 1}\phi^{\mathrm{inh}}|
\lesssim\left|Z^{\leq 1}\kl(|\psi_t|^p\kr)\right|
\lesssim
\kl(|\hpartial{u} Z^{\leq 1}\psi|+|\hpartial{v} Z^{\leq 1}\psi|+|\psi_r|\kr)
|\psi_t|^{p-1},
\end{Eq*}
\cref{claim:initial-energy-bound} and \eqref{eq:bootstrap-assumption},
together with the pointwise estimate \eqref{eq:radial-pointwise-estimate}
on $\psi$, imply
\begin{Eq*}
E_2[Z^{\leq 1}\phi^{\mathrm{inh}}](T)&
\lesssim\varepsilon^{2p}M^{2p}+\varepsilon^{p-1}M^{p-1}
\sum\Low{\star\in\{\mathrm{ext},\mathrm{mid},\mathrm{in}\}\\k=1,2,3,4,5,6}
I_{\phi;k}^\star,\\
I_{\phi;k}^\star:=&
\int_{S_0^{T;\star}}
r^{-\frac{n-2}{2}(p-1)}\kl<u\kr>^{1-p}\kl<v\kr>^{-\frac12(p-1)}I_{\phi;k}',
\end{Eq*}
where
\begin{Eq}\label{eq:source-pairings-and-regions}
I_{\phi;1}':=&\kl<u\kr>|\hpartial{u} Z^{\leq 1}\phi^{\mathrm{inh}}||\hpartial{u} Z^{\leq 1}\psi|,
\quad
I_{\phi;2}':=\kl<u\kr>|\hpartial{u} Z^{\leq 1}\phi^{\mathrm{inh}}||\hpartial{v} Z^{\leq 1}\psi|,
\quad
I_{\phi;3}':=\kl<u\kr>|\hpartial{u} Z^{\leq 1}\phi^{\mathrm{inh}}||\psi_r|,\\
I_{\phi;4}':=&\kl<v\kr>|\hpartial{v} Z^{\leq 1}\phi^{\mathrm{inh}}||\hpartial{u} Z^{\leq 1}\psi|,
\quad
I_{\phi;5}':=\kl<v\kr>|\hpartial{v} Z^{\leq 1}\phi^{\mathrm{inh}}||\hpartial{v} Z^{\leq 1}\psi|,
\quad
I_{\phi;6}':=\kl<v\kr>|\hpartial{v} Z^{\leq 1}\phi^{\mathrm{inh}}||\psi_r|,\\
S_0^{T;\mathrm{ext}}:=&S_0^T\cap\{r>2(t+1)\},\quad
S_0^{T;\mathrm{mid}}:=S_0^T\cap
\kl\{\frac{t+1}{2}<r\leq2(t+1)\kr\},\\
S_0^{T;\mathrm{in}}:=&S_0^T\cap
\kl\{r\leq\frac{t+1}{2}\kr\}.
\end{Eq}

We first consider the far exterior region, where $r\approx\kl<v\kr>\approx\kl<u\kr>\approx\kl<r\kr>\gtrsim\kl<t\kr>$. 
Since $\kl<u\kr>\approx\kl<v\kr>\approx\kl<r\kr>$ on this region,
the time-slice term of $E_1[f](t)$ in
\eqref{eq:energy-definitions}
controls the $\kl<r\kr>$-weighted $L^2$ norms of
$\hpartial{u}f$ and $\hpartial{v}f$. For every $k=1,2,4,5$, noting that $E_2[Z^{\leq 1}\psi](T)\lesssim\varepsilon^2M^2$, we have
\begin{Eq*}
I_{\phi;k}^{\mathrm{ext}}
\lesssim&\int_{S_0^{T;\mathrm{ext}}}
\kl<r\kr>\kl<t\kr>^{-\frac{n+1}{2}(p-1)}
|\hpartial{u,v}Z^{\leq 1}\phi^{\mathrm{inh}}|
|\hpartial{u,v}Z^{\leq 1}\psi|\\
\lesssim&
\int_0^T
    \kl<t\kr>^{-\frac{n+1}{2}(p-1)}
    \kl(
        \int_{\{r>2(t+1)\}}
            \kl<r\kr>|\hpartial{u,v}Z^{\leq 1}\phi^{\mathrm{inh}}|^2\kr
    )^{\frac12}
    \cdot
    \kl(
        \int_{\{r>2(t+1)\}}
            \kl<r\kr>|\hpartial{u,v}Z^{\leq 1}\psi|^2
    \kr)^{\frac12}
\d t\\
\lesssim&
\Big(E_2[Z^{\leq 1}\phi^{\mathrm{inh}}](T)\Big)^{\frac12}
\Big(E_2[Z^{\leq 1}\psi](T)\Big)^{\frac12}
\int_0^T
\kl<t\kr>^{s_p-p}\d t\\
\lesssim&\varepsilon M A_{s_p-p+1}(T)\Big(E_2[Z^{\leq 1}\phi^{\mathrm{inh}}](T)\Big)^{\frac12}.
\end{Eq*}
For $k=3,6$, using \eqref{eq:radial-pointwise-estimate} on $\psi$, we have
\begin{Eq*}
I_{\phi;k}^{\mathrm{ext}}
\lesssim&
\varepsilon M
\int_0^T
\int_{\{r>2(t+1)\}}
r^{1-\frac{n+1}{2}p}
|\hpartial{u,v}Z^{\leq 1}\phi^{\mathrm{inh}}|
\d x\d t\\
\lesssim&
\varepsilon M
\int_0^T
\kl(
\int_{\{r>2(t+1)\}}
\kl<r\kr>|\hpartial{u,v}Z^{\leq 1}\phi^{\mathrm{inh}}|^2\d x
\kr)^{\frac12}\\
&\qquad\qquad\times
\kl(
\int_{2(t+1)}^\infty
r^{-(n+1)(p-1)-1}\d r
\kr)^{\frac12}\d t\\
\lesssim&\varepsilon M \Big(E_2[Z^{\leq 1}\phi^{\mathrm{inh}}](T)\Big)^{\frac12}
\int_0^T
\kl<t\kr>^{s_p-p}\d t\\
\lesssim&
\varepsilon MA_{s_p-p+1}(T) \Big(E_2[Z^{\leq 1}\phi^{\mathrm{inh}}](T)\Big)^{\frac12}.
\end{Eq*}

We next turn to the interior region, where
$\kl<u\kr>\approx\kl<v\kr>\approx\kl<t\kr>\gtrsim\kl<r\kr>$.
For $k=1,2,4,5$, we have
\begin{Eq*}
I_{\phi;k}^{\mathrm{in}}
\lesssim&\int_{S_0^{T;\mathrm{in}}}
r^{-\frac{n-2}{2}(p-1)}
\kl<r\kr>^{-\frac32(p-1)}
\kl(Q[Z^{\leq 1}\phi^{\mathrm{inh}}]Q[Z^{\leq 1}\psi]\kr)^{1/2}\\
\lesssim&
\kl(
\int_{S_0^{T;\mathrm{in}}}
r^{-\frac{n-2}{2}(p-1)}
\kl<r\kr>^{-\frac32(p-1)}
Q[Z^{\leq 1}\phi^{\mathrm{inh}}]
\kr)^{\frac12}\times
\kl(
\int_{S_0^{T;\mathrm{in}}}
r^{-\frac{n-2}{2}(p-1)}
\kl<r\kr>^{-\frac32(p-1)}
Q[Z^{\leq 1}\psi]
\kr)^{\frac12},
\end{Eq*}
where $Q[f]$ is defined by \eqref{eq:Q-definition}. Then, 
since $p<n/(n-2)$,
applying
\eqref{eq:E2-weighted-spacetime} with
\begin{Eq*}
\delta:=1-\frac{n-2}{2}(p-1)>0,
\qquad
\eta:=s_p-p+1,
\end{Eq*}
for $k=1,2,4,5$, we find
\begin{Eq*}
I_{\phi;k}^{\mathrm{in}}
\lesssim&
A_{s_p-p+1}(T)
\Big(E_2[Z^{\leq 1}\phi^{\mathrm{inh}}](T)\Big)^{\frac12}
\Big(E_2[Z^{\leq 1}\psi](T)\Big)^{\frac12}\\
\lesssim&\varepsilon M A_{s_p-p+1}(T)\Big(E_2[Z^{\leq 1}\phi^{\mathrm{inh}}](T)\Big)^{\frac12}.
\end{Eq*}
For $k=3,6$, using \eqref{eq:radial-pointwise-estimate} and the same choice of
$\delta$ and $\eta$, we have
\begin{Eq*}
I_{\phi;k}^{\mathrm{in}}
\lesssim&
\varepsilon M
\int_{S_0^{T;\mathrm{in}}}
r^{-\frac{n-2}{2}p}
\kl<t\kr>^{-\frac{3p-1}{2}}
\Big(Q[Z^{\leq 1}\phi^{\mathrm{inh}}]\Big)^{\frac12}\\
\lesssim&
\varepsilon M
A_{s_p-p+1}(T)^{\frac12}
\Big(E_2[Z^{\leq 1}\phi^{\mathrm{inh}}](T)\Big)^{\frac12}
\kl(
\int_0^T
\kl<t\kr>^{1-3p}
\int_0^{(t+1)/2}
r^{n-1-\frac{n-2}{2}(p+1)}
\kl<r\kr>^{\frac32(p-1)}
\d r\d t
\kr)^{\frac12}\\
\lesssim&
\varepsilon M
A_{s_p-p+1}(T)^{\frac12}
\Big(E_2[Z^{\leq 1}\phi^{\mathrm{inh}}](T)\Big)^{\frac12}
\kl(
\int_0^T
\kl<t\kr>^{s_p-p}\d t
\kr)^{\frac12}\\
\lesssim&
\varepsilon M A_{s_p-p+1}(T) \Big(E_2[Z^{\leq 1}\phi^{\mathrm{inh}}](T)\Big)^{\frac12}.
\end{Eq*}

Now, we consider the middle region, where
$r\approx\kl<r\kr>\approx\kl<v\kr>\gtrsim\kl<u\kr>$. 
For $I_{\phi;1}^{\mathrm{mid}}$, since $\kl<u\kr>\geq1$, we obtain
\begin{Eq*}
I_{\phi;1}^{\mathrm{mid}}
\lesssim&
\int_0^{(3T+2)/2}
\kl<v\kr>^{-\frac{n-1}{2}(p-1)}
\kl(
\int_{\rHb_{-v,\min\{v,T-v\}}^v}
\kl<u\kr>|\hpartial{u} Z^{\leq 1}\phi^{\mathrm{inh}}|^2
\kr)^{\frac12}\\
&\qquad\qquad\qquad\times
\kl(
\int_{\rHb_{-v,\min\{v,T-v\}}^v}
\kl<u\kr>|\hpartial{u} Z^{\leq 1}\psi|^2
\kr)^{\frac12}\d v\\
\lesssim&
\Big(E_2[Z^{\leq 1}\phi^{\mathrm{inh}}](T)\Big)^{\frac12}\Big(E_2[Z^{\leq 1}\psi](T)\Big)^{\frac12}\int_0^{(3T+2)/2}
\kl<v\kr>^{-\frac{n-1}{2}(p-1)}\d v\\
\lesssim&\varepsilon M A_{s_p}(T)\Big(E_2[Z^{\leq 1}\phi^{\mathrm{inh}}](T)\Big)^{\frac12}.
\end{Eq*}

For $I_{\phi;2}^{\mathrm{mid}}$, the Cauchy--Schwarz inequality gives
\begin{Eq*}
I_{\phi;2}^{\mathrm{mid}}
\lesssim&
\kl(
\int_{S_0^{T;\mathrm{mid}}}
\kl<u\kr>^{4-p}\kl<v\kr>^{-(n-1)(p-1)-1}
|\hpartial{u} Z^{\leq 1}\phi^{\mathrm{inh}}|^2
\kr)^{\frac12}\\
&\quad\times
\kl(
\int_{S_0^{T;\mathrm{mid}}}
\kl<u\kr>^{-p}
\kl<v\kr>
|\hpartial{v} Z^{\leq 1}\psi|^2
\kr)^{\frac12}\\
\equiv&I_{\phi;21}^{1/2}I_{\phi;22}^{1/2}.
\end{Eq*}
For $I_{\phi;21}$, since $\kl<u\kr>\lesssim\kl<v\kr>$ and $p<n/(n-2)\leq3$ for all $n\geq3$, we obtain
\begin{Eq*}
I_{\phi;21}
\lesssim&
\int_0^{(3T+2)/2}\kl<v\kr>^{2s_p-p}
\kl(
\int_{\rHb^v_{-v,\min\{v,T-v\}}}
\kl<u\kr>|\hpartial{u} Z^{\leq 1}\phi^{\mathrm{inh}}|^2
\kr)\d v\\
\lesssim&E_2[Z^{\leq 1}\phi^{\mathrm{inh}}](T)
\int_0^{(3T+2)/2}
\kl<v\kr>^{2s_p-p}\d v\\
\lesssim&
A_{2s_p-p+1}(T)E_2[Z^{\leq 1}\phi^{\mathrm{inh}}](T).
\end{Eq*}
On the other hand, we have
\begin{Eq*}
I_{\phi;22}
\lesssim&
E_2[Z^{\leq 1}\psi](T)
\int_{-(T+2)/2}^{T/2}
\kl<u\kr>^{-p}\d u
\lesssim \varepsilon^2 M^2.
\end{Eq*}
Hence,
\begin{Eq*}
I_{\phi;2}^{\mathrm{mid}}
\lesssim
\varepsilon M A_{2s_p-p+1}(T)^{\frac12}
\Big(E_2[Z^{\leq 1}\phi^{\mathrm{inh}}](T)\Big)^{\frac12}.
\end{Eq*}

For $I_{\phi;3}^{\mathrm{mid}}$, using
\eqref{eq:radial-pointwise-estimate} and noting that $1-2p<-1$, we obtain
\begin{Eq*}
I_{\phi;3}^{\mathrm{mid}}
\lesssim&
\varepsilon M
\int_{S_0^{T;\mathrm{mid}}}
\kl<u\kr>^{1-p}
\kl<v\kr>^{-\frac{n-1}{2}p}
|\hpartial{u} Z^{\leq 1}\phi^{\mathrm{inh}}|
\d x\d t\\
\lesssim&
\varepsilon M
\int_0^{(3T+2)/2}
\kl<v\kr>^{-\frac{n-1}{2}p}
\kl(
\int_{\rHb_{-v,\min\{v,T-v\}}^v}
\kl<u\kr>|\hpartial{u} Z^{\leq 1}\phi^{\mathrm{inh}}|^2
\kr)^{\frac12}\\
&\qquad\times
\kl(
\int_{-v}^{\min\{v,T-v\}}
\kl<u\kr>^{1-2p}\kl<v\kr>^{n-1}\d u
\kr)^{\frac12}\d v\\
\lesssim&
\varepsilon M \Big(E_2[Z^{\leq 1}\phi^{\mathrm{inh}}](T)\Big)^{\frac12}
\int_0^{(3T+2)/2}
\kl<v\kr>^{s_p-1}\d v\\
\lesssim&
\varepsilon M A_{s_p}(T)\Big(E_2[Z^{\leq 1}\phi^{\mathrm{inh}}](T)\Big)^{\frac12}.
\end{Eq*}

For $I_{\phi;4}^{\mathrm{mid}}$, we use the additional $u$-decay in the first
square factor and write
\begin{Eq*}
I_{\phi;4}^{\mathrm{mid}}
\lesssim&
\kl(
\int_{S_0^{T;\mathrm{mid}}}
\kl<u\kr>^{1-2p}\kl<v\kr>
|\hpartial{v} Z^{\leq 1}\phi^{\mathrm{inh}}|^2
\kr)^{\frac12}\\
&\times
\kl(
\int_{S_0^{T;\mathrm{mid}}}
\kl<u\kr>\kl<v\kr>^{-(n-1)(p-1)+1}
|\hpartial{u} Z^{\leq 1}\psi|^2
\kr)^{\frac12}\\
\equiv&I_{\phi;41}^{1/2}I_{\phi;42}^{1/2}.
\end{Eq*}
We similarly have
\begin{Eq*}
I_{\phi;41}
\lesssim&E_2[Z^{\leq 1}\phi^{\mathrm{inh}}](T)\int_{-(T+2)/2}^{T/2}\kl<u\kr>^{1-2p}\d u
\lesssim E_2[Z^{\leq 1}\phi^{\mathrm{inh}}](T).
\end{Eq*}
On the other hand, we obtain
\begin{Eq*}
I_{\phi;42}
\lesssim&
\int_0^{(3T+2)/2}
\kl<v\kr>^{-(n-1)(p-1)+1}
\kl(
\int_{\rHb^v_{-v,\min\{v,T-v\}}}
\kl<u\kr>|\hpartial{u} Z^{\leq 1}\psi|^2
\kr)\d v\\
\lesssim&E_2[Z^{\leq 1}\psi](T)
\int_0^{(3T+2)/2}
\kl<v\kr>^{2s_p-1}\d v\\
\lesssim&\varepsilon^2 M^2 A_{2s_p}(T).
\end{Eq*}
Therefore, $I_{\phi;4}^{\mathrm{mid}}
\lesssim
\varepsilon M A_{2s_p}(T)^{\frac12}\Big(E_2[Z^{\leq 1}\phi^{\mathrm{inh}}](T)\Big)^{\frac12}$.

Next, for $I_{\phi;5}^{\mathrm{mid}}$, using
$\kl<v\kr>\gtrsim\kl<u\kr>$, we obtain
\begin{Eq*}
I_{\phi;5}^{\mathrm{mid}}
\lesssim&
\int_{-(T+2)/2}^{T/2}
\kl<u\kr>^{-\frac{n+1}{2}(p-1)}
\kl(
\int_{\rH_u^{|u|,T-u}}
\kl<v\kr>|\hpartial{v} Z^{\leq 1}\phi^{\mathrm{inh}}|^2
\kr)^{\frac12}\\
&\qquad\qquad\qquad\times
\kl(
\int_{\rH_u^{|u|,T-u}}
\kl<v\kr>
|\hpartial{v} Z^{\leq 1}\psi|^2
\kr)^{\frac12}\d u\\
\lesssim&
\Big(E_2[Z^{\leq 1}\phi^{\mathrm{inh}}](T)\Big)^{\frac12}\Big(E_2[Z^{\leq 1}\psi](T)\Big)^{\frac12}
\int_{-(T+2)/2}^{T/2}\kl<u\kr>^{s_p-p}\d u\\
\lesssim&\varepsilon M A_{s_p-p+1}(T)
\Big(E_2[Z^{\leq 1}\phi^{\mathrm{inh}}](T)\Big)^{\frac12}.
\end{Eq*}

Finally, for $I_{\phi;6}^{\mathrm{mid}}$, using
\eqref{eq:radial-pointwise-estimate}, we obtain
\begin{Eq*}
I_{\phi;6}^{\mathrm{mid}}
\lesssim&
\varepsilon M
\int_{S_0^{T;\mathrm{mid}}}
\kl<u\kr>^{-p}
\kl<v\kr>^{1-\frac{n-1}{2}p}
|\hpartial{v} Z^{\leq 1}\phi^{\mathrm{inh}}|
\d x\d t\\
\lesssim&
\varepsilon M
\int_{-(T+2)/2}^{T/2}
\kl<u\kr>^{-p}
\kl(
\int_{\rH_u^{|u|,T-u}}
\kl<v\kr>|\hpartial{v} Z^{\leq 1}\phi^{\mathrm{inh}}|^2
\kr)^{\frac12}\\
&\qquad\times
\kl(
\int_{|u|}^{T-u}
\kl<v\kr>^{2s_p-1}\d v
\kr)^{\frac12}\d u\\
\lesssim&
\varepsilon M A_{2s_p}(T)^{\frac12}
\Big(E_2[Z^{\leq 1}\phi^{\mathrm{inh}}](T)\Big)^{\frac12}
\int_{-(T+2)/2}^{T/2}
\kl<u\kr>^{-p}\d u\\
\lesssim&
\varepsilon M A_{2s_p}(T)^{\frac12}\Big(E_2[Z^{\leq 1}\phi^{\mathrm{inh}}](T)\Big)^{\frac12}.
\end{Eq*}

For fixed $\eta_1\leq\eta_2$, the definition
\eqref{eq:A-eta-definition} gives
\begin{Eq*}
A_{\eta_1}(T)\lesssim A_{\eta_2}(T),
\qquad
A_{2\eta}(T)^{\frac12}\lesssim A_\eta(T),
\end{Eq*}
where the implicit constants may depend on the fixed exponents. Combining the
estimates in the three regions gives
\begin{Eq}\label{eq:phi-inhomogeneous-energy-inequality}
E_2[Z^{\leq 1}\phi^{\mathrm{inh}}](T)
\lesssim\varepsilon^{2p}M^{2p}+\varepsilon^p M^p A_{s_p}(T)\Big(E_2[Z^{\leq 1}\phi^{\mathrm{inh}}](T)\Big)^{\frac12},
\end{Eq}
which leads to
\begin{Eq}\label{eq:phi-inhomogeneous-energy-bound}
E_2[Z^{\leq 1}\phi^{\mathrm{inh}}](T)
\lesssim\varepsilon^{2p}M^{2p}A_{s_p}(T)^2.
\end{Eq}

\subsection{Control of
\texorpdfstring{$E_2[Z^{\leq 1}\psi^{\mathrm{inh}}]$}{E2[Z<=1 psi inh]}
and closure of the bootstrap}

We now turn to the $\psi$-component. Applying
\eqref{eq:E2-estimate} to $Z^{\leq 1}\psi^{\mathrm{inh}}$, as in the preceding
subsection, we obtain
\begin{Eq}\label{eq:psi-inhomogeneous-energy}
E_2[Z^{\leq 1}\psi^{\mathrm{inh}}](T)
\lesssim&
E_0[Z^{\leq 1}\psi^{\mathrm{inh}}]
+\int_{S_0^T}
|\square Z^{\leq 1}\psi^{\mathrm{inh}}|
\kl(
\kl<u\kr>|\hpartial{u} Z^{\leq 1}\psi^{\mathrm{inh}}|
+\kl<v\kr>|\hpartial{v} Z^{\leq 1}\psi^{\mathrm{inh}}|
\kr).
\end{Eq}
Moreover,
\begin{Eq*}
|\square Z^{\leq 1}\psi^{\mathrm{inh}}|
\lesssim&
\left|Z^{\leq1}\kl(|\phi_t|^q\kr)\right|
\lesssim
\kl(
|\hpartial{u} Z^{\leq 1}\phi|
+|\hpartial{v} Z^{\leq 1}\phi|
+|\phi_r|
\kr)
|\phi_t|^{q-1}.
\end{Eq*}

At this point, we use
\begin{Eq*}
\phi=\phi^{\mathrm{hom}}+\phi^{\mathrm{inh}},
\qquad
|\phi_t|^{q-1}
\lesssim|\phi_t^{\mathrm{hom}}|^{q-1}+|\phi_t^{\mathrm{inh}}|^{q-1}.
\end{Eq*}
Using these relations in the preceding source estimate yields four types
of contributions: hom--hom, hom--inh, inh--hom, and inh--inh.
For simplicity of exposition, we only consider the two endpoint cases.
The mixed contributions can be controlled by the same regionwise
H\"older and Cauchy--Schwarz estimates, since the homogeneous and
inhomogeneous components satisfy weighted estimates with the same spatial
weights, and Young's inequality reduces the resulting coefficient
products to the two endpoint cases. Thus, no additional smallness
condition is required. We leave the details to the interested reader.

The purely homogeneous contribution is estimated as in the preceding
subsection, with $(p,\psi,\phi^{\mathrm{inh}})$ replaced by
$(q,\phi^{\mathrm{hom}},\psi^{\mathrm{inh}})$, since
$\phi^{\mathrm{hom}}$ satisfies the same uniform weighted energy bounds
as $\psi$ in the preceding subsection and hence enjoys the same
pointwise decay estimates.
More precisely, the estimates leading to
\eqref{eq:phi-inhomogeneous-energy-inequality} show that the purely
homogeneous contribution to the spacetime integral in
\eqref{eq:psi-inhomogeneous-energy}
is bounded by a constant multiple of
\begin{Eq*}
&\varepsilon^qM^q
A_{1-\frac{n-1}{2}(q-1)}(T)
\Big(E_2[Z^{\leq1}\psi^{\mathrm{inh}}](T)\Big)^{\frac12}\lesssim
\varepsilon^{q+1}M^{q+1}
A_{1-\frac{n-1}{2}(q-1)}(T),
\end{Eq*}
where the last inequality follows from
\eqref{eq:bootstrap-assumption}.
We next consider the fully inhomogeneous contribution, writing
\begin{Eq*}
|\square Z^{\leq 1}\psi^{\mathrm{inh}}|
\lesssim&
\kl(
|\hpartial{u} Z^{\leq 1}\phi^{\mathrm{inh}}|
+|\hpartial{v} Z^{\leq 1}\phi^{\mathrm{inh}}|
+|\phi_r^{\mathrm{inh}}|
\kr)
|\phi_t^{\mathrm{inh}}|^{q-1}\\
&+\text{terms containing at least one homogeneous factor}.
\end{Eq*}
Applying \eqref{eq:radial-pointwise-estimate} to
$\phi^{\mathrm{inh}}$ and using
\eqref{eq:phi-inhomogeneous-energy-bound}, we have
\begin{Eq*}
|\phi_r^{\mathrm{inh}}(u,v)|+|\phi_t^{\mathrm{inh}}(u,v)|&\lesssim
\varepsilon^p M^p r^{-\frac{n-2}{2}}A_{s_p}(t)
\kl<u\kr>^{-1}\kl<v\kr>^{-\frac12}.
\end{Eq*}
Arguing as in the preceding subsection, we therefore obtain
\begin{Eq*}
E_2[Z^{\leq 1}\psi^{\mathrm{inh}}](T)
\lesssim{}&\varepsilon^{2q}M^{2q}+\varepsilon^{p(q-1)}M^{p(q-1)}
\sum\Low{\star\in\{\mathrm{ext},\mathrm{mid},\mathrm{in}\}\\k=1,2,3,4,5,6}
I_{\psi;k}^\star\\
&+\text{terms containing at least one homogeneous factor},\\
I_{\psi;k}^\star:=&
\int_{S_0^{T;\star}}
r^{-\frac{n-2}{2}(q-1)}A_{s_p}(t)^{q-1}
\kl<u\kr>^{1-q}\kl<v\kr>^{-\frac12(q-1)}I_{\psi;k}',
\end{Eq*}
where the regions $S_0^{T;\star}$ are defined in
\eqref{eq:source-pairings-and-regions} and
\begin{Eq*}
I_{\psi;1}':=&\kl<u\kr>|\hpartial{u} Z^{\leq 1}\psi^{\mathrm{inh}}||\hpartial{u} Z^{\leq 1}\phi^{\mathrm{inh}}|,
\quad
I_{\psi;2}':=\kl<u\kr>|\hpartial{u} Z^{\leq 1}\psi^{\mathrm{inh}}||\hpartial{v} Z^{\leq 1}\phi^{\mathrm{inh}}|,
\quad
I_{\psi;3}':=\kl<u\kr>|\hpartial{u} Z^{\leq 1}\psi^{\mathrm{inh}}||\phi_r^{\mathrm{inh}}|,\\
I_{\psi;4}':=&\kl<v\kr>|\hpartial{v} Z^{\leq 1}\psi^{\mathrm{inh}}||\hpartial{u} Z^{\leq 1}\phi^{\mathrm{inh}}|,
\quad
I_{\psi;5}':=\kl<v\kr>|\hpartial{v} Z^{\leq 1}\psi^{\mathrm{inh}}||\hpartial{v} Z^{\leq 1}\phi^{\mathrm{inh}}|,
\quad
I_{\psi;6}':=\kl<v\kr>|\hpartial{v} Z^{\leq 1}\psi^{\mathrm{inh}}||\phi_r^{\mathrm{inh}}|.
\end{Eq*}

Inspection of the estimates in the preceding subsection shows that the
exterior and interior regions in fact yield stronger bounds. In the
middle region, $I_{\psi;1}^{\mathrm{mid}}$ and
$I_{\psi;3}^{\mathrm{mid}}$ attain the borderline growth, while
$I_{\psi;4}^{\mathrm{mid}}$ and $I_{\psi;6}^{\mathrm{mid}}$ yield
better bounds when $s_p=0$, and
$I_{\psi;2}^{\mathrm{mid}}$ and $I_{\psi;5}^{\mathrm{mid}}$ are always
better. Therefore, for simplicity of exposition, we only present the
estimates for $I_{\psi;1}^{\mathrm{mid}}$,
$I_{\psi;3}^{\mathrm{mid}}$, and $I_{\psi;4}^{\mathrm{mid}}$ below.
These estimates contain all the borderline cases and suffice to obtain
the desired sharp bound for $E_2[Z^{\leq 1}\psi^{\mathrm{inh}}](T)$.
The remaining terms are estimated in the same way and yield no worse time
growth.

In the middle region we have $r\approx\kl<r\kr>\approx\kl<v\kr>\gtrsim\kl<u\kr>\gtrsim 1$.
For $I_{\psi;1}^{\mathrm{mid}}$, arguing as for
$I_{\phi;1}^{\mathrm{mid}}$ in the preceding subsection, we obtain
\begin{Eq*}
I_{\psi;1}^{\mathrm{mid}}
\lesssim&
\int_0^{(3T+2)/2}
A_{s_p}(v)^{q-1}
\kl<v\kr>^{-\frac{n-1}{2}(q-1)}
\kl(
\int_{\rHb_{-v,\min\{v,T-v\}}^v}
\kl<u\kr>
|\hpartial{u} Z^{\leq 1}\psi^{\mathrm{inh}}|^2
\kr)^{\frac12}\\
&\qquad\qquad\qquad\times
\kl(
\int_{\rHb_{-v,\min\{v,T-v\}}^v}
\kl<u\kr>
|\hpartial{u} Z^{\leq 1}\phi^{\mathrm{inh}}|^2
\kr)^{\frac12}\d v.
\end{Eq*}
By \eqref{eq:bootstrap-assumption},
\begin{Eq*}
\kl(
\int_{\rHb_{-v,\min\{v,T-v\}}^v}
\kl<u\kr>
|\hpartial{u} Z^{\leq 1}\psi^{\mathrm{inh}}|^2
\kr)^{\frac12}
\lesssim
\varepsilon M.
\end{Eq*}
On the other hand, since $t=u+v\leq\min\{2v,T\}$ on the above
null segment, \eqref{eq:phi-inhomogeneous-energy-bound} gives
\begin{Eq*}
\kl(
\int_{\rHb_{-v,\min\{v,T-v\}}^v}
\kl<u\kr>
|\hpartial{u} Z^{\leq 1}\phi^{\mathrm{inh}}|^2
\kr)^{\frac12}
&\lesssim
\Big(E_2[Z^{\leq 1}\phi^{\mathrm{inh}}](\min\{2v,T\})\Big)^{\frac12}
\lesssim
\varepsilon^pM^pA_{s_p}(v).
\end{Eq*}
Therefore,
\begin{Eq*}
I_{\psi;1}^{\mathrm{mid}}
\lesssim
\varepsilon^{p+1}M^{p+1}\widetilde I_{\psi;1}^{\mathrm{mid}},\qquad
\widetilde I_{\psi;1}^{\mathrm{mid}}:=
\int_0^{(3T+2)/2}
A_{s_p}(v)^q
\kl<v\kr>^{-\frac{n-1}{2}(q-1)}
\d v.
\end{Eq*}
When $p>\frac{n+1}{n-1}$ so that $A_{s_p}(v)=1$,
$q\geq p$ implies $\Gamma(p,q)<0$ and $-\frac{n-1}{2}(q-1)<-1$, and hence
$\widetilde I_{\psi;1}^{\mathrm{mid}}\lesssim 1$.
When $p=\frac{n+1}{n-1}$ so that $A_{s_p}(v)=\ln(2+v)$,  we have
\begin{Eq*}
\widetilde I_{\psi;1}^{\mathrm{mid}}
\lesssim&
\int_0^{(3T+2)/2}
\ln^q(2+v)\kl<v\kr>^{-\frac{n-1}{2}(q-1)}\d v\\
\lesssim&
\begin{cases}
\ln^{p+1}(2+T),&q=p~\Rightarrow~ \Gamma(p,q)=0,\\
1,&q>p~\Rightarrow~ \Gamma(p,q)<0.
\end{cases}
\end{Eq*}
Finally, when $p<\frac{n+1}{n-1}$ so that $A_{s_p}(v)=\kl<v\kr>^{s_p}$,  we have
\begin{Eq*}
\widetilde I_{\psi;1}^{\mathrm{mid}}
\lesssim&
\int_0^{(3T+2)/2}\kl<v\kr>^{qs_p-\frac{n-1}{2}(q-1)}\d v
=\int_0^{(3T+2)/2}
\kl<v\kr>^{(pq-1)\Gamma(p,q)-1}\d v\\
\lesssim&
\begin{cases}
1,
&\Gamma(p,q)<0,\\
\ln(2+T),
&\Gamma(p,q)=0~\Rightarrow p<q,\\
\kl<T\kr>^{(pq-1)\Gamma(p,q)},
&\Gamma(p,q)>0.
\end{cases}
\end{Eq*}
Combining these cases, we conclude that $I_{\psi;1}^{\mathrm{mid}}\lesssim\varepsilon^{p+1}M^{p+1}\Lambda(T)$, where
\begin{Eq*}
\Lambda(T):=
\begin{cases}
1,
&\Gamma(p,q)<0,\\
\ln(2+T),
&\Gamma(p,q)=0,\quad p<q,\\
\ln^{p+1}(2+T),
&\Gamma(p,q)=0,\quad p=q,\\
\kl<T\kr>^{(pq-1)\Gamma(p,q)},
&\Gamma(p,q)>0.
\end{cases}
\end{Eq*}

Next, for $I_{\psi;3}^{\mathrm{mid}}$, using
\eqref{eq:radial-pointwise-estimate} for $\phi^{\mathrm{inh}}$, we obtain
\begin{Eq*}
I_{\psi;3}^{\mathrm{mid}}
\lesssim&
\varepsilon^pM^p
\int_{S_0^{T;\mathrm{mid}}}
A_{s_p}(t)^q
\kl<u\kr>^{1-q}
\kl<v\kr>^{-\frac{n-1}{2}q}
|\hpartial{u} Z^{\leq 1}\psi^{\mathrm{inh}}|
\d x\d t\\
\lesssim&
\varepsilon^pM^p
\int_0^{(3T+2)/2}
A_{s_p}(v)^q
\kl<v\kr>^{-\frac{n-1}{2}q}
\kl(
\int_{\rHb_{-v,\min\{v,T-v\}}^v}
\kl<u\kr>
|\hpartial{u} Z^{\leq 1}\psi^{\mathrm{inh}}|^2
\kr)^{\frac12}\\
&\qquad\times
\kl(
\int_{-v}^{\min\{v,T-v\}}
\kl<u\kr>^{1-2q}
\kl<v\kr>^{n-1}
\d u
\kr)^{\frac12}
\d v.
\end{Eq*}
Since $q>1$, we have
\begin{Eq*}
\int_{-v}^{\min\{v,T-v\}}
\kl<u\kr>^{1-2q}\d u
\lesssim1.
\end{Eq*}
Together with \eqref{eq:bootstrap-assumption}, this gives
\begin{Eq*}
I_{\psi;3}^{\mathrm{mid}}
&\lesssim
\varepsilon^{p+1}M^{p+1}
\int_0^{(3T+2)/2}
A_{s_p}(v)^q
\kl<v\kr>^{-\frac{n-1}{2}(q-1)}
\d v\\
&=
\varepsilon^{p+1}M^{p+1}
\widetilde I_{\psi;1}^{\mathrm{mid}}
\lesssim
\varepsilon^{p+1}M^{p+1}\Lambda(T).
\end{Eq*}

For $I_{\psi;4}^{\mathrm{mid}}$, as in the estimate of
$I_{\phi;4}^{\mathrm{mid}}$ in the preceding subsection, we write
\begin{Eq*}
I_{\psi;4}^{\mathrm{mid}}
\lesssim&
\kl(
\int_{S_0^{T;\mathrm{mid}}}
\kl<u\kr>^{1-2q}\kl<v\kr>
|\hpartial{v} Z^{\leq 1}\psi^{\mathrm{inh}}|^2
\kr)^{\frac12}\\
&\times
\kl(
\int_{S_0^{T;\mathrm{mid}}}
\kl<u\kr>
A_{s_p}(t)^{2(q-1)}
\kl<v\kr>^{-(n-1)(q-1)+1}
|\hpartial{u} Z^{\leq 1}\phi^{\mathrm{inh}}|^2
\kr)^{\frac12}\\
\equiv&
I_{\psi;41}^{\frac12}I_{\psi;42}^{\frac12}.
\end{Eq*}
Since $1-2q<-1$, the bootstrap assumption gives
\begin{Eq*}
I_{\psi;41}
&\lesssim
E_2[Z^{\leq 1}\psi^{\mathrm{inh}}](T)
\int_{-(T+2)/2}^{T/2}
\kl<u\kr>^{1-2q}\d u\lesssim
\varepsilon^2M^2.
\end{Eq*}
On the other hand, using $A_{s_p}(t)\lesssim A_{s_p}(v)$ in the
middle region and arguing as in the estimate of
$I_{\psi;1}^{\mathrm{mid}}$, we obtain
\begin{Eq*}
I_{\psi;42}
\lesssim&
\int_0^{(3T+2)/2}
A_{s_p}(v)^{2(q-1)}
\kl<v\kr>^{-(n-1)(q-1)+1}\\
&\qquad\qquad\times
\kl(
\int_{\rHb_{-v,\min\{v,T-v\}}^v}
\kl<u\kr>
|\hpartial{u} Z^{\leq 1}\phi^{\mathrm{inh}}|^2
\kr)\d v\\
\lesssim&
\varepsilon^{2p}M^{2p}
\int_0^{(3T+2)/2}
A_{s_p}(v)^{2q}
\kl<v\kr>^{-(n-1)(q-1)+1}\d v.
\end{Eq*}
Consequently, $I_{\psi;4}^{\mathrm{mid}}
\lesssim
\varepsilon^{p+1}M^{p+1}
\widetilde I_{\psi;4}^{\mathrm{mid}}$,
where
\begin{Eq*}
\widetilde I_{\psi;4}^{\mathrm{mid}}
:=
\kl(
\int_0^{(3T+2)/2}
A_{s_p}(v)^{2q}
\kl<v\kr>^{-(n-1)(q-1)+1}\d v
\kr)^{\frac12}.
\end{Eq*}
We next estimate $\widetilde I_{\psi;4}^{\mathrm{mid}}$.
When $p>\frac{n+1}{n-1}$ so that $A_{s_p}(v)=1$,
$q\geq p$ implies $\Gamma(p,q)<0$ and $-(n-1)(q-1)+1<-1$,
and hence $\widetilde I_{\psi;4}^{\mathrm{mid}}
\lesssim1$.
When $p=\frac{n+1}{n-1}$ so that
$A_{s_p}(v)=\ln(2+v)$, we have
\begin{Eq*}
\widetilde I_{\psi;4}^{\mathrm{mid}}
\lesssim&
\kl(
\int_0^{(3T+2)/2}
\ln^{2q}(2+v)
\kl<v\kr>^{-(n-1)(q-1)+1}
\d v
\kr)^{\frac12}\\
\lesssim&
\begin{cases}
\ln^{p+\frac12}(2+T),
&q=p~\Rightarrow~\Gamma(p,q)=0,\\
1,
&q>p~\Rightarrow~\Gamma(p,q)<0.
\end{cases}
\end{Eq*}
Finally, when $p<\frac{n+1}{n-1}$ so that
$A_{s_p}(v)=\kl<v\kr>^{s_p}$, we have
\begin{Eq*}
\widetilde I_{\psi;4}^{\mathrm{mid}}
\lesssim&
\kl(
\int_0^{(3T+2)/2}
\kl<v\kr>^{2qs_p-(n-1)(q-1)+1}
\d v
\kr)^{\frac12}\\
=&
\kl(
\int_0^{(3T+2)/2}
\kl<v\kr>^{2(pq-1)\Gamma(p,q)-1}
\d v
\kr)^{\frac12}\\
\lesssim&
\begin{cases}
1,
&\Gamma(p,q)<0,\\
\ln^{\frac12}(2+T),
&\Gamma(p,q)=0~\Rightarrow~p<q,\\
\kl<T\kr>^{(pq-1)\Gamma(p,q)},
&\Gamma(p,q)>0.
\end{cases}
\end{Eq*}
Combining these cases, we conclude that
$\widetilde I_{\psi;4}^{\mathrm{mid}}
\lesssim\Lambda(T)$,
and hence $
I_{\psi;4}^{\mathrm{mid}}
\lesssim
\varepsilon^{p+1}M^{p+1}\Lambda(T)$.

Combining the displayed estimates with \cref{claim:initial-energy-bound},
we obtain, for some constant $C>0$,
\begin{Eq*}
E_0[Z^{\leq 1}\psi^{\mathrm{inh}}]
&\leq C\varepsilon^{2q}M^{2q},\\
\varepsilon^{p(q-1)}M^{p(q-1)}
\sum\Low{\star\in\{\mathrm{ext},\mathrm{mid},\mathrm{in}\}\\k=1,2,3,4,5,6}
I_{\psi;k}^\star
&\leq C\varepsilon^{pq+1}M^{pq+1}\Lambda(T).
\end{Eq*}
Including the homogeneous and mixed contributions estimated above gives
\begin{Eq*}
E_2[Z^{\leq 1}\psi^{\mathrm{inh}}](T)
\lesssim{}&\varepsilon^{2q}M^{2q}
+\varepsilon^{q+1}M^{q+1}
A_{1-\frac{n-1}{2}(q-1)}(T)+\varepsilon^{pq+1}M^{pq+1}\Lambda(T)\\
\lesssim{}&\varepsilon^2M^2\kl(
\varepsilon^{2q-2}M^{2q-2}
+\kl(\varepsilon^{pq-1}M^{pq-1}\Lambda(T)\kr)^{\frac{q-1}{pq-1}}
+\varepsilon^{pq-1}M^{pq-1}\Lambda(T)
\kr),
\end{Eq*}
where we used the fact that
$A_{1-\frac{n-1}{2}(q-1)}(T)
\lesssim\Lambda(T)^{\frac{q-1}{pq-1}}$ for $p\leq q$.
Hence these contributions impose no additional smallness condition, and it
suffices to verify
\begin{Eq}\label{eq:bootstrap-smallness-conditions}
\varepsilon^{2q-2}M^{2q-2}\ll1,
\qquad
\varepsilon^{pq-1}M^{pq-1}\Lambda(T)\ll1.
\end{Eq}
Choose $0<\varepsilon_0\leq1$ sufficiently small, depending only on $n,p,q$
and the fixed constant $M$, so that all the required time-independent
smallness conditions hold. We decrease $\varepsilon_0$ and the lifespan
constant $c$ further below whenever necessary.

If $\Gamma(p,q)<0$, then $\Lambda(T)=1$. Hence, for
$0<\varepsilon\leq\varepsilon_0$,
the conditions in \eqref{eq:bootstrap-smallness-conditions} hold without any
restriction on $T$.

Suppose next that $\Gamma(p,q)=0$ and $p\neq q$. We have
$\Lambda(T)=\ln(2+T)$. For
\begin{Eq*}
T\leq\exp\kl(c\varepsilon^{-(pq-1)}\kr),
\end{Eq*}
we have
\begin{Eq*}
\ln(2+T)\leq\ln3+c\varepsilon^{-(pq-1)}.
\end{Eq*}
After decreasing $\varepsilon_0$ and choosing $c>0$ sufficiently small,
depending only on $n,p,q$ and the fixed constants, the conditions in
\eqref{eq:bootstrap-smallness-conditions} follow.

Suppose now that $\Gamma(p,q)=0$ and $p=q$. Since
$pq-1=(p-1)(p+1)$ and $\Lambda(T)=\ln^{p+1}(2+T)$, for
\begin{Eq*}
T\leq\exp\kl(c\varepsilon^{-(p-1)}\kr),
\end{Eq*}
we have
\begin{Eq*}
\ln(2+T)\leq\ln3+c\varepsilon^{-(p-1)}.
\end{Eq*}
After decreasing $\varepsilon_0$ and choosing $c>0$ sufficiently small, the
conditions in \eqref{eq:bootstrap-smallness-conditions} again follow.

Finally, suppose that $\Gamma(p,q)>0$. We have
$\Lambda(T)=\kl<T\kr>^{(pq-1)\Gamma(p,q)}$. After decreasing
$\varepsilon_0$ further if necessary, the conditions in
\eqref{eq:bootstrap-smallness-conditions} hold for
\begin{Eq*}
T\leq c\varepsilon^{-\frac{1}{\Gamma(p,q)}},
\end{Eq*}
for a sufficiently small constant $c>0$.

In all four cases, the preceding estimates give
\begin{Eq*}
E_2[Z^{\leq 1}\psi^{\mathrm{inh}}](t)\leq2\varepsilon^2M^2,
\qquad 0\leq t\leq T,
\end{Eq*}
for every $T<T_b$. This strictly improves
\eqref{eq:bootstrap-assumption}, so weighted energy continuity on compact
subintervals of $[0,T_*)$ gives $T_b=\min\{T_*,T_\varepsilon\}$.
If $T_*<T_\varepsilon$, the preceding energy bounds control the unweighted
derivative norms, while the $L^2$ norms of $\phi$ and $\psi$ remain bounded
by time integration of $\phi_t$ and $\psi_t$. Thus the Sobolev norms in
\cref{prop:local-well-posedness} remain uniformly bounded for $0\leq t<T_*$.
The continuation criterion in that proposition then extends the solution
beyond $T_*$, a contradiction. Therefore $T_*\geq T_\varepsilon$.
Finally, finite propagation and the uniform weighted a priori bounds allow
us to pass from compactly supported data to general data in
$\mathcal X^{\mathrm{rad}}$ by the standard approximation argument.
This proves all the cases in \cref{thm:main} and completes the proof.

\subsection{Initial energy bound}

\begin{proof}[Proof of \cref{claim:initial-energy-bound}]
We first consider the homogeneous parts. By definition,
\begin{Eq*}
\phi^{\mathrm{hom}}(0)&=\varepsilon f_\phi,
&
\phi_t^{\mathrm{hom}}(0)&=\varepsilon g_\phi,\\
\psi^{\mathrm{hom}}(0)&=\varepsilon f_\psi,
&
\psi_t^{\mathrm{hom}}(0)&=\varepsilon g_\psi.
\end{Eq*}
Since $\square\phi^{\mathrm{hom}}=\square\psi^{\mathrm{hom}}=0$, we also have
\begin{Eq*}
\phi_{tt}^{\mathrm{hom}}(0)
=&\varepsilon\Delta f_\phi,
\qquad
\psi_{tt}^{\mathrm{hom}}(0)
=\varepsilon\Delta f_\psi.
\end{Eq*}
Moreover, on the initial hypersurface,
\begin{Eq*}
S f^{\mathrm{hom}}(0)
=&r\partial_r f^{\mathrm{hom}}(0),\\
\partial_t(S f^{\mathrm{hom}})(0)
=&f_t^{\mathrm{hom}}(0)
+r\partial_r f_t^{\mathrm{hom}}(0),
\end{Eq*}
while
\begin{Eq*}
L_i f^{\mathrm{hom}}(0)
=&x_i f_t^{\mathrm{hom}}(0),\\
\partial_t(L_i f^{\mathrm{hom}})(0)
=&\partial_i f^{\mathrm{hom}}(0)
+x_i f_{tt}^{\mathrm{hom}}(0).
\end{Eq*}
Therefore, by the definition of $E_0$ in \eqref{eq:energy-definitions},
\begin{Eq*}
\Big(E_0[Z^{\leq 1}\phi^{\mathrm{hom}}]\Big)^{\frac12}
+\Big(E_0[Z^{\leq 1}\psi^{\mathrm{hom}}]\Big)^{\frac12}
\lesssim
\varepsilon
\left\|
\kl<r\kr>^{\frac32}
\kl(
\nabla^{\leq2}f_\phi,
\nabla^{\leq1}g_\phi,
\nabla^{\leq2}f_\psi,
\nabla^{\leq1}g_\psi
\kr)
\right\|_{L^2}
\lesssim
\varepsilon\|\boldsymbol f\|_{\mathcal X}.
\end{Eq*}

We next consider the inhomogeneous parts. Set
\begin{Eq*}
N_\phi:=|\psi_t|^p,
\qquad
N_\psi:=|\phi_t|^q.
\end{Eq*}
Since
\begin{Eq*}
\square\phi^{\mathrm{inh}}=N_\phi,
\qquad
\square\psi^{\mathrm{inh}}=N_\psi,
\end{Eq*}
and both inhomogeneous parts have zero Cauchy data, we have
\begin{Eq*}
\phi^{\mathrm{inh}}(0)
=\phi_t^{\mathrm{inh}}(0)
=&0,
&
\phi_{tt}^{\mathrm{inh}}(0)
=&\varepsilon^p|g_\psi|^p,\\
\psi^{\mathrm{inh}}(0)
=\psi_t^{\mathrm{inh}}(0)
=&0,
&
\psi_{tt}^{\mathrm{inh}}(0)
=&\varepsilon^q|g_\phi|^q.
\end{Eq*}
Consequently, the fields
\begin{Eq*}
f^{\mathrm{inh}},
\qquad
\partial_i f^{\mathrm{inh}},
\qquad
S f^{\mathrm{inh}}
\end{Eq*}
have zero Cauchy data. The only nonzero Cauchy data among the first-order
commuted fields are
\begin{Eq*}
\partial_t f^{\mathrm{inh}}(0)&=0,
&
\partial_t^2 f^{\mathrm{inh}}(0)&=N_f(0),\\
L_i f^{\mathrm{inh}}(0)&=0,
&
\partial_t(L_i f^{\mathrm{inh}})(0)&=x_iN_f(0),
\end{Eq*}
where
\begin{Eq*}
N_\phi(0)=\varepsilon^p|g_\psi|^p,
\qquad
N_\psi(0)=\varepsilon^q|g_\phi|^q.
\end{Eq*}
It follows that
\begin{Eq*}
\Big(E_0[Z^{\leq 1}\phi^{\mathrm{inh}}]\Big)^{\frac12}
&\lesssim
\varepsilon^p
\left\|
\kl<r\kr>^{\frac32}|g_\psi|^p
\right\|_{L^2},\\
\Big(E_0[Z^{\leq 1}\psi^{\mathrm{inh}}]\Big)^{\frac12}
&\lesssim
\varepsilon^q
\left\|
\kl<r\kr>^{\frac32}|g_\phi|^q
\right\|_{L^2}.
\end{Eq*}

For every $1<s<n/(n-2)$, the Sobolev inequality applied to
$\kl<r\kr>^{3/(2s)}h$ gives
\begin{Eq*}
\left\|
\kl<r\kr>^{\frac{3}{2s}}h
\right\|_{L^{2s}}
\lesssim
\left\|
\nabla^{\leq1}
\kl(
\kl<r\kr>^{\frac{3}{2s}}h
\kr)
\right\|_{L^2}
\lesssim
\left\|
\kl<r\kr>^{\frac32}
\nabla^{\leq1}h
\right\|_{L^2}.
\end{Eq*}
Hence,
\begin{Eq*}
\left\|
\kl<r\kr>^{\frac32}|g_\psi|^p
\right\|_{L^2}
=&
\left\|
\kl<r\kr>^{\frac{3}{2p}}g_\psi
\right\|_{L^{2p}}^p
\lesssim
\left\|
\kl<r\kr>^{\frac32}
\nabla^{\leq1}g_\psi
\right\|_{L^2}^p,\\
\left\|
\kl<r\kr>^{\frac32}|g_\phi|^q
\right\|_{L^2}
=&
\left\|
\kl<r\kr>^{\frac{3}{2q}}g_\phi
\right\|_{L^{2q}}^q
\lesssim
\left\|
\kl<r\kr>^{\frac32}
\nabla^{\leq1}g_\phi
\right\|_{L^2}^q.
\end{Eq*}
Therefore,
\begin{Eq*}
\Big(E_0[Z^{\leq 1}\phi^{\mathrm{inh}}]\Big)^{\frac12}
&\lesssim
\varepsilon^p
\|\boldsymbol f\|_{\mathcal X}^p,\\
\Big(E_0[Z^{\leq 1}\psi^{\mathrm{inh}}]\Big)^{\frac12}
&\lesssim
\varepsilon^q
\|\boldsymbol f\|_{\mathcal X}^q.
\end{Eq*}

Using the normalization
\begin{Eq*}
\|\boldsymbol f\|_{\mathcal X}\leq1,
\end{Eq*}
we conclude that
\begin{Eq*}
\Big(E_0[Z^{\leq 1}\phi^{\mathrm{hom}}]\Big)^{\frac12}
+\Big(E_0[Z^{\leq 1}\psi^{\mathrm{hom}}]\Big)^{\frac12}
&\lesssim\varepsilon,\\
\Big(E_0[Z^{\leq 1}\phi^{\mathrm{inh}}]\Big)^{\frac12}
&\lesssim\varepsilon^p,\\
\Big(E_0[Z^{\leq 1}\psi^{\mathrm{inh}}]\Big)^{\frac12}
&\lesssim\varepsilon^q.
\end{Eq*}
Using the above decompositions and $\varepsilon\leq1$, we also obtain
\begin{Eq*}
\Big(E_0[Z^{\leq 1}\phi]\Big)^{\frac12}+\Big(E_0[Z^{\leq 1}\psi]\Big)^{\frac12}
\lesssim\varepsilon.
\end{Eq*}
This completes the proof.
\end{proof}

\subsection*{Acknowledgment}
The first author was supported by NSFC 12501292.
The second author was supported by NSFC 12101558.
The third author was supported by NSFC 12671293.

\bibliography{bib}

\end{document}